\documentclass[11pt,a4paper]{article}
\usepackage{amsfonts}
\usepackage{amsmath,amssymb}
\usepackage{mathrsfs}
\usepackage{hyperref}
\usepackage{graphicx}
\usepackage{float}
\usepackage{amsmath,amssymb,amsthm,amsfonts,xcolor,graphicx,units,xfrac}
\usepackage{geometry}

\newtheorem{thm}{Theorem}[section]
\newtheorem{cor}[thm]{Corollary}
\newtheorem{lem}[thm]{Lemma}
\newtheorem{prop}[thm]{Proposition}

\numberwithin{equation}{section}\allowdisplaybreaks
\def\leq{\leqslant}

\def\leq{\leqslant}
\def\geq{\geqslant}

\def\Real{{\mathbb{R}}}

\def\F{{\mathscr{F}}}

\begin{document}

\title{\large\bf Global Well-Posedness and Global Attractor for Two-dimensional Zakharov-Kuznetsov Equation}
\author{\normalsize \bf  Minjie Shan \footnote{Corresponding author. The author was supported by CSC, grant 201606010025} \\
\footnotesize
\it School of Mathematical Sciences, Peking University, Beijing 100871, P.R. China\\
\footnotesize
\texttt{ Emails: smj@pku.edu.cn}
} \maketitle

\thispagestyle{empty}
\begin{abstract}
  The initial value problem for two-dimensional Zakharov-Kuznetsov equation is shown to be globally well-posed in $H^{s}({\Real^2})$ for all $\frac{5}{7}<s<1$ via using $I$-method in the context of atomic spaces. By means of the increment of modified energy, the exsitence of global attractor for
  weakly damped, forced Zakharov-Kuznetsov equation is also established in $H^{s}({\Real^2})$ for $\frac{10}{11}<s<1$. \\


\end{abstract}

\section{Introduction and Main Results}

We consider the Cauchy problem for the symmetrized two-dimensional Zakharov-Kuznetsov (ZK) equation
\begin{equation}
  \left\{
   \begin{aligned}
   &u_{t} +( \partial^3_{x}+\partial^3_{y}) u +( \partial_{x}+\partial_{y}) u^2 = 0, \ \ \ (x,y)\in\Real^2, \ t\geq0     \quad \\
   &u(x,y,0)=u_0(x,y)\in H^s(\Real^2), \label{ZK} \\
   \end{aligned}
   \right.
  \end{equation}
where  $u=u(x,y,t)$ is a real-valued function.

The ZK equation was initially decuded as a model of nonlinear unidirectional ion-acoustic wave propagation in a magnetized plasma by Zakharov
and Kuznetsov \cite{ZK74}. It may be treated as a higher dimensional generalization of the
Korteweg-de Vries (KdV) equation. For more details we refer to the papers \cite{EWLKHS82, LLS13} about the two-dimensional ZK equation appearing here in physical
circumstances. 

Eventhough the ZK equation is not completely integrable, there still exists two conserved quantities for
the flow of ZK,

\begin{align}
M(u)(t)= \int_{\Real^2} u^2(x,y,t) dxdy = \int_{\Real^2} u^2_0(x,y) dxdy = M(u_0) \label{M()}
\end{align}

and
\begin{align}
E(u)(t)= \int_{\Real^2} \frac{1}{2} | \nabla u|^2 -\frac{1}{2} \partial_{x}u\partial_{y}u - \frac{1}{3} u^3 dxdy = E(u_0).\label{E()}
\end{align}

Faminskii \cite{F95} firstly obtained the local well-posedness for the two-dimensional ZK equation in the energy space $H^1(\Real^2)$ by making use of local smoothing effects together with a maximal function estimate for the linearized equation. This method was inspired by the one given by Kenig, Ponce and Vega \cite{KPV93} when dealing with the local well-posedness for KdV equation. With the help of the  $L^2$ and $H^1$ conservation
laws, he proved the global well-posedness for ZK equation additionally. Following this idea, Linares and Pastor \cite{LP09} optimized the proof of Faminskii to show local well-posedness in $H^s(\Real^2)$ for $ s > \frac{3}{4}$. Gr\"unrock and Herr \cite{GH14} along with Molinet and Pilod \cite{MP15} proved the local well-posedness in a larger data spaces $H^s(\Real^2)$ for $ s > \frac{1}{2}$  by taking advantage of the Fourier restriction norm method and one kind of sharp Strichartz estimates. On the basis of the method of Gr\"unrock and Herr, in \cite{Shan1} we recently improved the local well-posedness result to $ B^\frac{1}{2}_{2,1}(\Real^2)$ via applying frequency decomposion as well as atomic spaces introduced by Koch and Tataru. Actually, it is L. Molinet and D. Pilod who originality applied this crucial tool of atomic spaces to ZK equation. They obtained the global well-posedness for ZK equation in $H^s(\Real^3)$ for  $s>1$ by using the conservation laws, doubling time, the argument of Bourgain (see \cite{Bourgain 2a}) and atomic spaces. What's more, as the other part of  the same paper \cite{Shan1} we utilized $I$-method to prove global well-posedness in $H^s(\Real^2)$ for $  \frac{11}{13}<s<1$.

As mentioned above, ZK equation is an asymptotic description of the propagation of nonlinear ion-acoustic wave in a magnetized plasma. However, one cannot totally neglect external excitation and energy dissipation mechanism in reality (see \cite{A1984}, \cite{OS1970} and \cite{Wu1987}). Therefore, we
would like to consider the following Zakharov-Kuznetsov equation with weak dissipation and forcing terms
\begin{equation}
  \left\{
   \begin{aligned}
   &u_{t} +( \partial^3_{x}+\partial^3_{y}) u +( \partial_{x}+\partial_{y}) u^2+\gamma u =f, \ \ \ (x,y)\in\Real^2, \  t\geq0     \quad \\
   &u(x,y,0)=u_0(x,y)\in H^s(\Real^2), \label{dfZK} \\
   \end{aligned}
   \right.
  \end{equation}
where the damping parameter $\gamma$ is a positive constant and the external forcing term $f\in H^1(\Real^2)$ is independent of $t$.

Global attractor is a bounded subset which is invariant by the flow and attracts all trajectories when $t$ approaches to $+\infty$ under the corresponding topology. Moise and Rosa studied the regularity of the global attractor of a weakly damped, forced KdV equation in \cite{MR1997}. While Haraux \cite{Haraux1985} and Goubet \cite{Goubet1996} proved the asymptotic smoothing effect for dissipative wave equation and nonlinear Schr\"{o}dinger eqution respectively. As we know that the conservation laws will become invalid in low regularity situation, Colliander, Keel, Staffilani, Takaoka and Tao \cite{CKSTT02,CKSTT03} developed a theory of almost conserved quantities starting with energy in order to extend the local solutions globally in time for dispersive equations. Applying $I$-method to weakly damped, forced KdV equation, Tsugawa \cite{Tsugawa} proved the existence of global attractor on Sobolev spaces of negative index. Prashant \cite{Pra} lately showed an analogue below energy spaces for mKdV equation.

In this article, we refine the global well-posedness for ZK equation through pushing our previous result to $  \frac{5}{7}<s$. The new ingredient here is a closer investigation into how to use the symmetrization and atomic space to acquire a better control of the modified energy. Whereafter the increment of the modified energy helps us gain the existence of global weak attractor below  $H^1$ for weakly damped, forced ZK equation. The advantage of $U^2$, $V^2$ spaces is to obtain sharp estimates in the time variable, which makes the proof simpler than when we use $X^{s, 1/2+}$. For instance, as regards local well-posedness for  ZK equation in $ H^{\frac{1}{2}}(\Real^2)$,  to some extent it corresponds to the non-admissible endpoint Strichartz estimates. However this can be compensated partly by the decay brought from utilizing atomic space. Especially, the use of $U^2$, $V^2$ spaces will assist us to study the problem on the 2D torus or the $\mathbb{R} \times \mathbb{T}$, which may be nice as a future problem (see also  \cite{MP15} for  ZK equation on $\mathbb{R} \times \mathbb{T}$).

We now state the main results of this paper.

\begin{thm}
The initial value problem \eqref{ZK} is globally well-posed in $H^s(\Real^2)$ for $\frac{5}{7}<s<1$.
\end{thm}
 \emph{Remark} $1.$ The increment in $E(I_Nu)$ on the right hand of \eqref {E incresement} is  $N^{-\frac{1}{2}}$. One could repeat the almost
 conservation law argument to prove global well-posedness of \eqref{ZK} for all  $ \frac{1}{2}<s<1$ upon gaining $N^{-1}$.

\begin{thm} \label{Global Attactor}
  Let $\frac{10}{11}<s<1$. Then, there exists a semigroup $A(t)$ and maps $M_1$ and $M_2$ such that $A(t)u_0$ is the unique solution of \eqref{dfZK} satisfying
\begin{align}
  A(t)u_0=M_1(t)u_0+M_2(t)u_0 \label{U(t)} \\
 \sup \limits_{t>T_1} \|M_1(t)u_0\|_{H^1} < K \label{M_1(t)}
\end{align}
and for $t>T_1$
\begin{align}
  \|M_2(t)u_0\|_{H^s} < K e^{-\gamma(t-T_1)} \label{M_2(t)}
\end{align}
where $T_1$ depends on $\|u_0\|_{H^s}$, $\|f\|_{H^1}$ and $\gamma$, the constant $K$ depends only on $\|f\|_{H^1}$ and $\gamma$.  
\end{thm}

From this theorem, we know that $M_1(t)$ is a bounded mapping and $M_2(t)$ converges uniformly to 0 in $H^s$. It means that the
semigroup $A(t)$ is asymptotically compact in the sense of weak topology. Therefore we gain the existence of the global attractor in $H^s$ from Theorem 1.1.1 in \cite{Temam}.

\begin{cor}
The global weak attractor for \eqref{dfZK} exisits in $H^s$ for  $\frac{10}{11}<s<1$.
\end{cor}

Note that
\begin{align}
  \partial_t\int u^2& =-2\int u( \partial^3_{x}+\partial^3_{y}) u-2\int u( \partial_{x}+\partial_{y}) u^2-2\gamma \int u^2 +2\int fu\nonumber \\
  &=-2\gamma \int u^2 +2\int fu \nonumber
\end{align}
Setting $h(t)=e^{2\gamma t}\|u\|^2_{L^2}$ and using Cauchy-Schwarz inequality in the last term, one can obtain
$$h'(t)\leqslant2e^{\gamma t}\|f\|_{L^2}\sqrt{h(t)}.$$
This implies that
$$\partial_t\sqrt{h(t)}\leqslant e^{\gamma t}\|f\|_{L^2}.$$
Hence, we have
\begin{align}
  \|u\|_{L^2}\leqslant e^{-\gamma t}\|u_0\|_{L^2}+\gamma^{-1}\|f\|_{L^2}\label{uL^2}.
\end{align}

\emph{Remark} $2.$ For the modified energy, we expect a similar exponential decay as \eqref{uL^2} (see Proposition \ref {Prop priori estimate2}).

\emph{Remark} $3.$ Yang \cite{Yang} upgrated Tsugawa's result through lowering $-\frac{3}{8}$ to $-\frac{1}{2}$ which can be seen as the critical
Sobolev index for KdV. However, for ZK equation and  mKdV equation it is not easy to find good cancellation of resonant parts as Schr\"{o}dinger equation (see \cite{CKSTT08}). Concerning the existence of global attrator, there still are much work to do for these two kinds of equations.
\\
\\\textbf{Organization of the paper.} In Section 2, we recall some useful propositions and estimates about $U^p$ and $V^p$. Then we prove global well-posedness by using $I$-method in atomic space in Section 3. Finally, Section 4 is devoted to a prior estimate and the proof of Theorem 1.2.

 We now list notations used throughout this paper.
 Let $c < 1$, $C>3$ denote universal constants. The notation $c+$ stands for $c+\epsilon$  some  $0<\epsilon\ll 1$. Similarly, we shall write $ c-= c-\epsilon$. We put $\langle a \rangle = (1+a^2)^ \frac{1}{2}$ for $a\in\Real$ and fix a smooth cut-off function
 $\chi\in C^\infty_0([-2,2])$ satisfing $\chi$ is even, nonnegative,  and $\chi=1$ on $[-1,1]$. We denote spatial variables by $x,y$ and their dual Fourier variables by $\xi,\eta$. As usual, $\tau$ is the dual variable of the time $t$. Let $\tilde f$ denote the Fourier transform of $f$ in both time and spatial variables. Let $\widehat f$ denote its Fourier transform only in space or in time. For $s\in\Real$, $I^s_x $ and $I^s_y $ denote the one-dimensional Riesz-potential operators of order $-s$ with respect to spatial variable $x$ and $y$ . We also write $\zeta=(\xi,\eta)$, $\lambda=(\xi,\eta,\tau)$ and $\mu =\tau -\xi^3-\eta^3 $ for brevity. We will make frequent use of the capital letters $N, M, N_1, N_2$ and $N_3$ which denote dyadic numbers and we  write $\sum_{N\geq 1}a_N=\sum_{n\in \mathbb{N}}a_{2^n}$ , $\sum_{N\geq M}a_N=\sum_{n\in \mathbb{N};2^n\geq M}a_{2^n}$ for dyadic summations.

\section{Function spaces and estimates}

In this section we introduce some properties of $U^p$ and $V^p$ spaces (see \cite{HaHeKo09,KoTa05,KoTa07,KoTa12}) which is another powful tool to handle low regularity well-posedness for dispersive equations.

Let $1\leq p <\infty$ and $\mathcal{Z}$ be the set of finite partitions $-\infty= t_0 <t_1<...< t_{K-1} < t_K =\infty$.

For any $\{t_k\}^K_{k=0} \subset \mathcal{Z}$ and $\{\phi_k\}^{K-1}_{k=0} \in L^2$ with $\sum^{K-1}_{k=0} \|\phi_k\|^p_2=1$, $\phi_0=0$. We call the function $a: \mathbb{R}\to L^2$ given by
$$
a= \sum^{K}_{k=1} 1_{[t_{k-1}, t_k)} \phi_{k-1}
$$
 a $U^p$-atom.  The atomic space is
$$
U^p=\left\{u= \sum^\infty_{j=1} \lambda_j a_j : a_j \ U^p-atom, \lambda_j \in \mathbb{C},  \sum^\infty_{j=1} |\lambda_j|<\infty  \right\}
$$
endowed with the norm
$$
\|u\|_{U^p}= \mathop{inf} \left\{\sum^\infty_{j=1} |\lambda_j| : u= \sum^\infty_{j=1} \lambda_j a_j ,a_j \ U^p-atom, \lambda_j \in \mathbb{C} \right\}.$$

$V^p$ is the normed space of all functions $v: \mathbb{R} \to L^2$ such that $\mathop{lim}\limits_{t\to \pm \infty} v(t)$ exist and for which the norm
$$
\|v\|_{V^p} = \mathop{sup}\limits_{\{t_k\}^K_{k=0} \in \mathcal{Z}} \left( \sum^K_{k=1} \|v(t_k)-v(t_{k-1})\|^p_{L^2}\right)^{\frac{1}{p}}  \label{v norm}
$$
is finite, where we use the convention that $v(-\infty) = \mathop{lim}\limits_{t\to - \infty} v(t)$ and $v(\infty)=0$. We denote $v \in V^p_-$ when $v(-\cdot) \in V^p $. Moreover, we define the closed subspace $V^p_{rc}$ $(V^p_{-,rc})$ of all right continuous functions in $V^p$ $(V^p_-)$.

The unitary operator $e^{tS}:L^2 \to L^2$ is defined to be the Fourier multiplier
$$\widehat {e^{tS}u_0}(\xi,\eta)=e^{it(\xi^3+\eta^3)}\widehat {u_0}(\xi,\eta),
$$
where $S=-\partial^3_x-\partial^3_y$. Let us define $U^p_S = e^{\cdot S} U^p$ with norm  $ \|u\|_{U^p_S} = \|e^{-\cdot S} u \|_{U^p}$ and $V^p_S = e^{\cdot S} V^p$ with norm  $ \|v\|_{V^p_S} = \|e^{-\cdot S} v \|_{V^p}$.

Given the Littlewood-Paley multipliers by 

$$\widehat{P_1 u}=\chi(2|\zeta|)\widehat{u}$$
 and $$\widehat{P_Nu}=\psi_N(|\zeta|)\widehat{u} \ \ \ for \ \  N\geq 2 \ ,$$
where $\psi(x) =\chi(x)-\chi(2x)$ and $\psi_N =\psi(N^{-1}\cdot)$, the smooth projections according to the dispersive relationship are defined by

$$\F (Q_M u)(\xi,\eta,\tau)=\psi_M(\tau)\tilde u(\xi,\eta,\tau),
$$
$$\F (Q^S_M u)(\xi,\eta,\tau)=\psi_M(\tau-\xi^3-\eta^3)\tilde u(\xi,\eta,\tau),
$$
as well as $Q^S_{\geq M}=\sum_{N\geq M}Q^S_N$ and $Q^S_{< M}=I-Q^S_{\geq M}$. Note that 
$Q^S_M=e^{\cdot S}Q_Me^{-\cdot S}.$\vspace{3mm}

Let's recall some useful results in $U^p$ and $V^p$.
\begin{prop}  \label{UV prop1} Let $1< p <q <\infty$ and $\frac{1}{p}+\frac{1}{p'}=1$. We have
\begin{itemize}
 \item[\rm (i)] $U^p$, $V^p$, $V^p_{rc}$, $V^p_{-}$ and $V^p_{-,rc}$ are Banach spaces,
\item[\rm (ii)] $U^p\subset V^p_{-,rc} \subset U^q$ ,
\item[\rm (iii)]$\|u\|_{U^p}=\sup \limits_{\|v\|_{V^{p'}}=1} |\int \langle u'(t), v(t)\rangle dt|$  if $u\in V^1_{-}(\subset U^p)$ is absolutely continuous on compact interval.
\end{itemize}
\end{prop}
\begin{lem}\label{UV lema1}
 We have
\begin{align}
  \|Q^S_{\geq M}u \|_{L^2(\Real^3) }  \lesssim M^{-\frac{1}{2}} \|u\|_{V^2_S}. \label{UV estimate 1}
\end{align}
\begin{align}
  \|Q^S_{\geq M}u \|_{U^p_S} \lesssim \|u\|_{U^p_S} , \ \ \  \|Q^S_{< M}u \|_{U^p_S} \lesssim \|u\|_{U^p_S}. \label{UV estimate 2}
\end{align}
\begin{align}
 \|Q^S_{\geq M}u \|_{V^p_S} \lesssim \|u\|_{V^p_S} , \ \ \  \|Q^S_{< M}u \|_{V^p_S} \lesssim \|u\|_{V^p_S}.\label{UV estimate 3}
\end{align}
\end{lem}\vspace{3mm}
Similarly to  Lemma 2.3 in \cite{GTV97} and Lemma 5.3 in \cite{WaHaHuGu11},  the extension principle also holds true for $U^p_S$ spaces.
\begin{prop}  \label{UV prop2}
 Let $ T_0: L^2 \times \cdots \times L^2 \to L^1_{loc}(\Real^2; \mathbb{C}) $ be a n-linear operator. Assume that for some $1< p,q <\infty$
$$\| T_0(e^{\cdot S}\phi_1, \cdots, e^{\cdot S}\phi_n )\|_{L^p_t(\Real; L^q_{x,y}(\Real^2))} \lesssim \prod_{j=1}^n\|\phi_j\|_{L^2}.$$
Then, there exists $ T: U^p_S \times \cdots \times U^p_S \to L^p_t(\Real; L^q_{x,y}(\Real^2)) $ satisfying
$$\| T(u_1, \cdots, u_n )\|_{L^p_t(\Real; L^q_{x,y}(\Real^2))} \lesssim \prod_{j=1}^n\|u_j\|_{U^p_S},$$
such that $ T(u_1, \cdots, u_n )(t)(x,y)=T_0(u_1(t), \cdots, u_n(t) )(x,y) \ a.e.$.
\end{prop}

Next we present the Strichartz estimates and bilinear estimates (see Section 5 in \cite{Shan1}).
\begin{lem} \label{in lema}
Given $N_1\geq N_2$. Assume that $(q,r)$ satisfy $\frac{3}{q}+\frac{2}{r}=1$ and $q>3$. Let $I^s_{x,-}$ be the bilinear operator with symbol $|\xi_1-\xi_2|^s$, i.e.

$$\F_{x,y}(I^s_{x,-}(f_1,f_2))(\xi,\eta)=\int_{\zeta=\zeta_1+\zeta_2}|\xi_1-\xi_2|^s \prod^2_{j=1}\widehat{f_j}(\xi_j,\eta_j).$$
Then, we have
\begin{align}
 \|\chi(\frac{t}{T})u \|_{L^2(\Real^3)}\lesssim T^{\frac{1}{2}}\|u\|_{V^2_S}, \label{estimate21}
\end{align}
\begin{align}
  \|u\|_{L^q_tL^r_{xy}}  \lesssim \|u\|_{U^q_S},\label{estimate22}
\end{align}
\begin{align}
  \|I^{\frac{1}{8}}_xI^{\frac{1}{8}}_ye^{tS}u_0\|_{L^4(\Real^3)}  \lesssim \|u_0\|_{L^2(\Real^2)},\label{estimate23}
\end{align}
\begin{align}
  \|I^{\frac{1}{2}}_xI^{\frac{1}{2}}_{x,-}(P_{N_1}e^{tS}u_0,P_{N_2}e^{tS}v_0) \|_{L^2(\Real^3)} \lesssim N_2^{\frac{1}{2}}\|P_{N_1}u_0\|_{L^2(\Real^2)}\|P_{N_2}v_0\|_{L^2(\Real^2)} , \label{estimate24}
\end{align}
Moreover, \eqref{estimate24} is equally valid with $x$ replaced by $y$.
\end{lem}
\begin{lem} \label{derivertive lema}
Let $0 \leq \epsilon < \frac{1}{2}$ and  $0 \leq \theta \leq 1 $. Then, 
\begin{align}
  \|\tilde D^{\frac{\theta\epsilon}{2}}e^{tS}u_0\|_{L^q_tL^r_{xy}}  \lesssim \|u_0\|_{L^2_{xy}},\label{estimate25}
\end{align}
where $\widehat{\tilde D^su}=|\xi+\eta|^s\hat u(\xi,\eta)$, $q=\frac{6}{\theta(2+\epsilon)}$ and $r=\frac{2}{1-\theta}$.
\end{lem}
{\bf Proof.}  Set $\tilde I_t(x,y)=\int_{\Real^2}|\xi+\eta|^{\epsilon+i\beta}e^{i[t(\xi^3+\eta^3)+(x\xi+y\eta)]}d\xi d\eta$. Let $ {\xi}'=4^{-\frac{1}{3}}(\xi+\eta)$ and $ {\eta}'=\sqrt{3}4^{-\frac{1}{3}}(\xi-\eta)$,  then we get
$$\tilde I_t(x,y)=cI_t(2^{-\frac{1}{3}}(x+y),2^{-\frac{1}{3}}3^{-\frac{1}{2}}(x-y)),$$ where $I_t(x,y)$ is defined as in Lemma 2.1 of \cite{LP09}.
This implies
$$\tilde I_t(x,y)\lesssim |t|^{-\frac{2+\epsilon}{3}}.$$
Hence, \eqref{estimate25} is a direct consequence of interpolation theorem and standard Stein-Thomas argument.

We can avail ourselves of the technique in \cite{MV15} to decompose the time cut-off into low- and high-frequency parts.

For any $\delta>0$, we write $1_\delta$ the characteristic function on $[0,\delta]$ and 
$$1_\delta=1^{low}_{\delta,\kappa}+1^{high}_{\delta,\kappa}, \ \ \widehat {1^{low}_{\delta,\kappa}}(\tau)=\chi(\tau/\kappa)\widehat {1_\delta}(\tau)$$
for some $\kappa >0$.
\begin{lem} \label{cutoff lema}
For any $\kappa, \delta>0$, it holds
\begin{align}
 \|1^{high}_{\delta,\kappa} \|_{L^{\frac{3}{2}}(\Real)}\lesssim \delta^{\frac{1}{3}}\kappa^{-\frac{1}{3}}, \label{estimate26}
\end{align}
\begin{align}
  \|1^{high}_{\delta,\kappa}\|_{L^{\infty}}\lesssim 1, \label{estimate27}
\end{align}
\begin{align}
  \|1^{low}_{\delta,\kappa}\|_{L^{\infty}}\lesssim 1. \label{estimate28}
\end{align}
\end{lem}

\section{Global well-posedness}
We introduce $I$-method. Given $m:\Real^{2k} \to \mathbb{C}$, $m$ is said to be symmetric if
$$m(\zeta_1, \cdots,\zeta_k)=m(\sigma(\zeta_1, \cdots,\zeta_k))$$
for all $\sigma \in S_k$, where $S_k$ is the permutation group for $k$ elements. The symbol $m$ can be symmetrized as following
$$[m]_{sym}(\zeta_1, \cdots,\zeta_k)=\frac{1}{k!}\sum_{\sigma \in S_k}m(\sigma(\zeta_1, \cdots,\zeta_k)).$$
We define a $k$-linear functional acting on $k$ functions $u_1, \cdots,u_k$ for each $m$
$$\Lambda_k(m;u_1, \cdots,u_k)=\int_{\zeta_1+\cdots+\zeta_k=0}m(\zeta_1, \cdots,\zeta_k)\widehat u_1(\zeta_k)\cdots\widehat u_k(\zeta_k).$$
Usually we write $\Lambda_k(m)$ instead of $\Lambda_k(m;u, \cdots,u)$ for convenience. Note that $\Lambda_k(m)= \Lambda_k([m]_{sym})$ from symmetries.

Given $s<1$, $N\gg1$ and a smooth, radially symmetric, nonincreasing function $m(\zeta)$ satisfying

$$m(\zeta)=
\left\{
\begin{array}{l}
1 \ \ \ \ \ \ \ \ \ \ \ \ |\zeta|\leq N \\
(\frac{|\zeta|}{N})^{s-1} \ \ \ \ |\zeta|\geq 2N
\end{array}
\right. ,$$
we define the Fourier multiplier operator
$$\widehat{If}(\zeta)=m(\zeta) \widehat f(\zeta).$$
Let $\lambda>0$ and $N'=\frac{N}{\lambda}$, we define
$$ m'(\zeta)=
\left\{
\begin{array}{l}
1 \ \ \ \ \ \ \ \ \ \ \ \ |\zeta|\leq N' \\
(\frac{|\zeta|}{N'})^{s-1} \ \ \ \ |\zeta|\geq 2N'
\end{array}
\right. $$
and rescaled operator
$$\widehat{I'f}(\zeta)=m'(\zeta) \widehat f(\zeta).$$

Notice that $$|\int \partial_{x}u\partial_{y}udxdy|\leq \frac{1}{2}\int |\nabla u|^2dxdy,$$ it's easy to prove  $\|u\|_{H^s}$ and $E(Iu)(t)$  are comparable as Proposition 3.2 in \cite{Shan1}.
\begin{prop} Let $\frac{5}{7}<s<1$, then
\begin{align}
|E(Iu)(t)|\lesssim N^{2(1-s)}\|u(t)\|^2_{\dot H^s(\Real^2)}+\|u(t)\|^3_{L^3(\Real^2)} \nonumber
\end{align}
\begin{align}
\|u(t)\|^2_{ H^s(\Real^2)}\lesssim |E(Iu)(t)|+\|u_0\|^2_{L^2(\Real^2)}+\|u_0\|^4_{L^2(\Real^2)}. \nonumber
\end{align}
\end{prop}

We start from local existence theroem. The work space is denoted by $Y^s $ which can be defined via the norm
$$\|u\|_{Y^s}=(\sum_N N^{2s}\|P_Nu\|^2_{U^2_S})^{\frac{1}{2}}.$$ Obviously, $Y^s \subset L^{\infty}H^s$ holds true. 
\begin{lem} \label{local lemma1}
  Let $N_1, N_2$ be dyadic numbers and $N_2\lesssim N_1$. Assume that $|\zeta_k|\sim N_k$ $(k=1,2)$ and $|\zeta_1+\zeta_2|\sim N_3$.
Denote $\Omega_1=\{\lambda_1+\lambda_2=\lambda,|\eta|\leq|\xi|,|\xi|\lesssim |\xi_1-\xi_2|\}$ and $\Omega_2=\{\lambda_1+\lambda_2=\lambda,|\eta|\leq|\xi|,|\xi|\sim |\xi_1|\sim|\xi_2|\lesssim|\eta_1|\sim|\eta_2|\}$. We have
\\$(i)$ if $N_1\sim N_3\gg N_2$, then
\begin{align}
 \|\int_{\Omega_1}(\xi+\eta)\frac{m(\zeta_1+\zeta_2)}{m(\zeta_1)m(\zeta_2)}\tilde u_{N_1}(\lambda_1)\tilde v_{N_2}(\lambda_2)d\lambda_1\|_{L^2_{\lambda}(\Real^3)}\lesssim N_2^{\frac{3}{2}-s}\|u_{N_1}\|_{U^2_S}\|v_{N_2}\|_{U^2_S}, \label{estimate31}
\end{align}
\\$(ii)$ if $N_1\sim  N_2\gtrsim N_3$, then
 \begin{align}
 \|\int_{\Omega_k}(\xi+\eta)\frac{m(\zeta_1+\zeta_2)}{m(\zeta_1)m(\zeta_2)}\tilde u_{N_1}(\lambda_1)\tilde v_{N_2}(\lambda_2)d\lambda_1\|_{L^2_{\lambda}(\Real^3)}\lesssim N_1^{\frac{1}{2}}C(N_1,N_3)\|u_{N_1}\|_{U^2_S}\|v_{N_2}\|_{U^2_S} \label{estimate32}
 \end{align}
 for $k=1,2$, where
 $$C(N_1,N_3)=\left\{\begin{array}{l}
    1 \ \ \ \ \ \ \ \ \ \ \ \ \ \ \ \ \ \  \ \ \ \ \ \ \ N\gg N_1\sim N_2\gg N_3 \\
 ({\frac{N_1}{N}})^{2-2s}\ \ \ \ \ \ \ \ \ \ \ \  \ \ \ \ N_1\sim N_2 \gtrsim N\gg N_3 \\
  
  N^{s-1}{N_1}^{2-2s}{N_3}^{s-1} \ \ \ \ N_1\sim N_2 \gg N_3\gtrsim N \\
  {N_1}^{1-s} \ \ \ \ \ \ \ \ \ \ \ \ \ \ \ \ \ \ \ N_1\sim N_2\sim N_3 
\end{array}
\right. .$$
\end{lem}
{\bf Proof.} It's suffices to show the estimates for free solutions by Proposition \ref{UV prop2}. To this end, we can prove it along the same line as the proof of Lemma 5.3 in \cite{Shan1}. 

If $N_1\sim N_3\gg N_2$, it is obvious that
$$|\frac{m(\zeta_1+\zeta_2)}{m(\zeta_1)m(\zeta_2)}|\lesssim 1\lor({\frac{N_2}{N}})^{1-s}\lesssim N_2^{1-s}.$$
Besides, we have $|\xi+\eta| \lesssim |\xi|\lesssim|\xi|^{\frac{1}{2}}|\xi_1-\xi_2|^{\frac{1}{2}}$ on $\Omega_1$. Therefore, \eqref{estimate24} gives
\begin{align}
& \|\int_{\Omega_1}(\xi+\eta)\frac{m(\zeta_1+\zeta_2)}{m(\zeta_1)m(\zeta_2)}\F (e^{tS}u_{0,N_1})(\lambda_1)) \F (e^{tS}v_{0,N_2})(\lambda_2)d\lambda_1\|_{L^2_{\lambda}(\Real^3)}\notag \\
   \lesssim &N_2^{1-s}\|\int_{\Omega_1}|\xi|^{\frac{1}{2}}|\xi_1-\xi_2|^{\frac{1}{2}}\F (e^{tS}u_{0,N_1})(\lambda_1)) \F (e^{tS}v_{0,N_2})(\lambda_2)d\lambda_1\|_{L^2_{\lambda}(\Real^3)}\notag \\
  \lesssim &N_2^{1-s}\|\ I^{\frac{1}{2}}_xI^{\frac{1}{2}}_{x,-}(e^{tS}u_{0,N_1},e^{tS}v_{0,N_2})\|_{L^2_{\lambda}(\Real^3)}\notag \\
  \lesssim & N_2^{\frac{3}{2}-s}\|u_{0,N_1}\|_{L^2_{x,y}}\|v_{0,N_2}\|_{L^2_{x,y}}.  \nonumber     
\end{align}
   
 If $N_1\sim  N_2\gtrsim N_3$, we can bound $\frac{m(\zeta_1+\zeta_2)}{m(\zeta_1)m(\zeta_2)}$ by $C(N_1,N_3)$.  Also, on hyperplane $\Omega_2$, $|\xi+\eta|\lesssim|\xi|\lesssim N_2^{\frac{1}{2}}|\xi_1|^{\frac{1}{8}}|\eta_1|^{\frac{1}{8}}|\xi_2|^{\frac{1}{8}}|\eta_2|^{\frac{1}{8}}$. From H\"{o}rder inequalities and
\eqref{estimate23}, we get
\begin{align}
& \|\int_{\Omega_2}(\xi+\eta)\frac{m(\zeta_1+\zeta_2)}{m(\zeta_1)m(\zeta_2)}\F (e^{tS}u_{0,N_1})(\lambda_1)\F (e^{tS}v_{0,N_2})(\lambda_2)d\lambda_1\|_{L^2_{\lambda}(\Real^3)}\notag \\
  \lesssim & N_1^{\frac{1}{2}}C(N_1,N_3)\|\int_{\Omega_2}\F({I^{\frac{1}{8}}_xI^{\frac{1}{8}}_ye^{tS}u_{0,N_1}})(\lambda_1)\F({I^{\frac{1}{8}}_xI^{\frac{1}{8}}_ye^{tS}v_{0,N_2}})(\lambda_2) d\lambda_1\|_{L^2_{\lambda}(\Real^3)}\notag \\
  \lesssim &  N_1^{\frac{1}{2}}C(N_1,N_3)\|I^{\frac{1}{8}}_xI^{\frac{1}{8}}_ye^{tS}u_{0,N_1}\|_{L^4(\Real^3)}\|I^{\frac{1}{8}}_xI^{\frac{1}{8}}_ye^{tS}v_{0,N_2}\|_{L^4(\Real^3)}\notag \\
  \lesssim &  N_1^{\frac{1}{2}}C(N_1,N_3)\|u_{0,N_1}\|_{L^2_{x,y}}\|v_{0,N_2}\|_{L^2_{x,y}}.  \nonumber     
\end{align}
Additionally, the estimate holds true on $\Omega_1$ in a similar way. Hence, the proof is completed.
    
\begin{lem} \label {local lemma2}  Let $N_1, N_2, N_3$ be dyadic numbers and $C(N_1,N_3)$ be given as in Lemma \ref{local lemma1}. Then, for all $0<\delta\leq1$ it holds
\\$(i)$ if $N_1\sim N_3\gg N_2$,
\begin{align}
&|\int_{\Real^3}\chi(\frac{t}{\delta})I(u_{N_1}v_{N_2})(\partial_x+\partial_y)w_{N_3}dxdydt|\notag \\
\lesssim& \delta^{\frac{1}{2}}N_2^{\frac{3}{2}-s}\|Iu_{N_1}\|_{U^2_S}\|Iv_{N_2}\|_{U^2_S}\|w_{N_3}\|_{V^2_S} \label{estimate33},\end{align}
\\$(ii)$  if $N_1\sim  N_2\gtrsim N_3$,
\begin{align}
&|\int_{\Real^3}\chi(\frac{t}{\delta})I(u_{N_1}v_{N_2})(\partial_x+\partial_y)w_{N_3}dxdydt|\notag \\
\lesssim& \delta^{\frac{1}{6}}N_1^{\frac{1}{2}}C(N_1,N_3)\|Iu_{N_1}\|_{U^2_S}\|Iv_{N_2}\|_{U^2_S}\|w_{N_3}\|_{V^2_S}.\label{estimate34}
\end{align}
\end{lem}
    {\bf Proof.} As the technique used in Proposition 5.4 of \cite{Shan1} is applicative here, we just give a sketch to avoid being lengthy and tedious. 
   
 Parseval formula shows that it suffices to control
$$\int_{\sum^3_{j=1}\lambda_j=0}(\xi_3+\eta_3)\frac{m(\zeta_1+\zeta_2)}{m(\zeta_1)m(\zeta_2)}\tilde u_{N_1}(\lambda_1)\tilde v_{N_2}(\lambda_2)\F (\chi_\delta w_{N_3})(\lambda_3)$$
by $\|u_{N_1}\|_{U^2_S}$, $\|v_{N_2}\|_{U^2_S}$ and $\|w_{N_3}\|_{V^2_S}$.

We assume that $|\eta_3|\leq|\xi_3|$ and denote the absolute value of the quatity above by $R$.

If $N_1\sim N_3\gg N_2$, then $|\xi_3|\lesssim |\xi_1-\xi_2|$.
Applying Cauchy-Schwarz inequality, \eqref {estimate31} and \eqref {estimate21} gives
\begin{align}
R&\lesssim \|\int_{\Omega_1}(\xi+\eta)\frac{m(\zeta_1+\zeta_2)}{m(\zeta_1)m(\zeta_2)}\tilde u_{N_1}(\lambda_1)\tilde v_{N_2}(\lambda_2)d\lambda_1\|_{L^2_{\lambda}(\Real^3)}\|\chi(\frac{t}{\delta})w_{N_3}\|_{L^2(\Real^3)}\notag \\
&\lesssim \delta^{\frac{1}{2}}N_2^{\frac{3}{2}-s}\|u_{N_1}\|_{U^2_S}\|v_{N_2}\|_{U^2_S}\|w_{N_3}\|_{V^2_S}. \nonumber
\end{align}

So as to prove \eqref{estimate34}, we need to split the domain of the integration into five regions $R=R_1+R_2+R_3+R_4+R_5$. One can assume $|\eta_1|\geq |\eta_2|$ by symmetry.
\\{\bf Region 1.} $|\xi_3|\lesssim |\xi_1-\xi_2|$

Using \eqref {estimate32}, we have
\begin{align}
R_1&\lesssim  \|\int_{\Omega_1}(\xi+\eta)\frac{m(\zeta_1+\zeta_2)}{m(\zeta_1)m(\zeta_2)}\tilde u_{N_1}(\lambda_1)\tilde v_{N_2}(\lambda_2)d\lambda_1\|_{L^2_{\lambda}(\Real^3)}\|\chi(\frac{t}{\delta})w_{N_3}\|_{L^2(\Real^3)}\notag \\
&\lesssim \delta^{\frac{1}{2}}N_1^{\frac{1}{2}}C(N_1,N_3)\|u_{N_1}\|_{U^2_S}\|v_{N_2}\|_{U^2_S}\|w_{N_3}\|_{V^2_S}.\nonumber
\end{align}
\\{\bf Region 2.} $|\xi_3|\gg |\xi_1-\xi_2|$ and $|\eta_1|\gg |\xi_3|$

This condition gives
$$|\eta_1|\sim |\eta_2|\gg |\xi_3|\sim|\xi_1|\sim|\xi_2|.$$
Applying Cauchy-Schwarz inequality, \eqref{estimate32} and \eqref{estimate21}, we obtain
\begin{align}
R_2&\lesssim  \|\int_{\Omega_2}(\xi+\eta)\frac{m(\zeta_1+\zeta_2)}{m(\zeta_1)m(\zeta_2)}\tilde u_{N_1}(\lambda_1)\tilde v_{N_2}(\lambda_2)d\lambda_1\|_{L^2(\Real^3)} \|\chi(\frac{t}{\delta})w_{N_3}\|_{L^2(\Real^3)}\notag \\
&\lesssim \delta^{\frac{1}{2}}N_1^{\frac{1}{2}}C(N_1,N_3)\|u_{N_1}\|_{U^2_S}\|v_{N_2}\|_{U^2_S}\|w_{N_3}\|_{V^2_S}.\nonumber
\end{align}
{\bf Region 3.} $|\xi_3|\gg |\xi_1-\xi_2|$ and $|\eta_1|\sim |\xi_3|\sim |\eta_2|$

We get the bound of $R_3$ like Region 2.
\\{\bf Region 4.} $|\xi_3|\gg |\xi_1-\xi_2|$ and $|\eta_1|\sim |\xi_3|\gg |\eta_2|$

We decompose $Id=Q^S_{<M}+Q^S_{\geq M}$ and divide the integral $R_4$ into eight pieces.
\\{\bf Case A.} $Q^S_j=Q^S_{<M}$ for $j=1,2,3$

We go a step further to decompose time into low- and high-frequency parts.

{\bf Case A(1).} All of these three are low-frequence

The integral is vanished.

{\bf Case A(2).} At least one of these three is high-frequence

We estimate when the first one is high-frequence
$$\int_{*}(\xi_3+\eta_3)\frac{m(\zeta_1+\zeta_2)}{m(\zeta_1)m(\zeta_2)}\F (1^{high}_{\delta,\kappa}Q^S_{<M}u_{N_1})\F (1_vQ^S_{<M}v_{N_2})\F (1_wQ^S_{<M}w_{N_3}),$$
where $1_v,1_w\in \left\{1^{high}_{\delta,\kappa},1^{low}_{\delta,\kappa}\right\}$.

H\"{o}lder inequalities, \eqref{estimate22}, Lemma \ref{cutoff lema} and \eqref {UV estimate 2} provide
\begin{align}
&|\int_{*}(\xi_3+\eta_3)\frac{m(\zeta_1+\zeta_2)}{m(\zeta_1)m(\zeta_2)}\F (1^{high}_{\delta,\kappa}Q^S_{<M}u_{N_1})\F (1_vQ^S_{<M}v_{N_2})\F (1_wQ^S_{<M}w_{N_3})| \notag \\  
\lesssim& \|1^{high}_{\delta,\kappa}Q^S_{<M}u^*_{N_1}\|_{L^{\frac{9}{7}}_tL^3_{xy}}\|1_vQ^S_{<M}v^*_{N_2}\|_{L^9_tL^3_{xy}}\|1_wQ^S_{<M}w^*_{N_3}\|_{L^9_tL^3_{xy}}    \notag \\
\lesssim& \|1^{high}_{\delta,L}\|_{L^{\frac{3}{2}}_t}\|Q^S_{<M}u^*_{N_1}\|_{L^9_tL^3_{xy}}\|1_v\|_{L^{\infty}_t}\|Q^S_{<M}v^*_{N_2}\|_{L^9_tL^3_{xy}}\|1_w\|_{L^{\infty}_t}
    \|Q^S_{<M}w^*_{N_3}\|_{L^9_tL^3_{xy}}   \notag \\
\lesssim& \delta^{\frac{1}{3}}N^{-\frac{1}{3}}_3N^{-\frac{2}{3}}_1\|Q^S_{<M}u^*_{N_1}\|_{U^9_S}\|Q^S_{<M}v^*_{N_2}\|_{U^9_S}\|Q^S_{<M}w^*_{N_3}\|_{U^9_S} \notag \\
\lesssim& \delta^{\frac{1}{3}}N^{-\frac{1}{3}}_3N^{-\frac{2}{3}}_1\|u^*_{N_1}\|_{U^2_S}\|v^*_{N_2}\|_{U^2_S}\|w^*_{N_3}\|_{V^2_S},\nonumber
\end{align}
where $u^*_{N_1}=\F^{-1}\frac{\widehat u(\zeta_1)}{m(\zeta_1)}$, $v^*_{N_2}=\F^{-1}\frac{\widehat v(\zeta_2)}{m(\zeta_2)}$ and $w^*_{N_3}=\F^{-1}(\xi_3+\eta_3)m(\zeta_3)\widehat w(\zeta_3)$.

From the definition of $U^2$ and $V^2$, there hold the facts
$$\|u^*_{N_1}\|_{U^2_S}\lesssim \frac{\|u_{N_1}\|_{U^2_S}}{m(N_1)},\ \  \|v^*_{N_2}\|_{U^2_S}\lesssim \frac{\|v_{N_2}\|_{U^2_S}}{m(N_2)}$$

and
$$\|w^*_{N_3}\|_{V^2_S}\lesssim N_3m(N_3)\|w_{N_3}\|_{V^2_S},$$
which of course help imply the desired up bound $$\delta^{\frac{1}{3}}N_1^{\frac{1}{2}}C(N_1,N_3)\|u_{N_1}\|_{U^2_S}\|v_{N_2}\|_{U^2_S}\|w_{N_3}\|_{V^2_S}.$$
{\bf Case B.} $Q^S_j=Q^S_{\geq M}$ for some $j=1,2,3$

We take $Q^S_3=Q^S_{\geq M}$ for instance to give the estimate.

Using H\"{o}lder inequalities, \eqref {UV estimate 1}, \eqref{estimate22} and \eqref {UV estimate 2}, we have
\begin{align}
&|\int_{*}(\xi_3+\eta_3)\frac{m(\zeta_1+\zeta_2)}{m(\zeta_1)m(\zeta_2)}\F (\chi_\delta Q^S_1u_{N_1})\F (\chi_\delta Q^S_2v_{N_2})
\F (\chi_\delta Q^S_{\geq M}w_{N_3})|\notag \\
\lesssim& \|\chi_\delta Q^S_1u^*_{N_1}\|_{L^4(\Real^3)}\|\chi_\delta Q^S_2v^*_{N_2}\|_{L^4(\Real^3)}\|Q^S_{\geq M}w^*_{N_3}\|_{L^2(\Real^3)} \notag \\
\lesssim& M^{-\frac{1}{2}}\|\chi_\delta\|^2_{L^{12}(\Real)}\|Q^S_1u^*_{N_1}\|_{L^6_tL^4_{x,y}}\|Q^S_2v^*_{N_2}\|_{L^6_tL^4_{x,y}}
\|w^*_{N_3}\|_{V^2_S}\notag \\
\lesssim& \delta^{\frac{1}{6}}N_1^{-1}N_3^{-\frac{1}{2}}\|Q^S_1u^*_{N_1}\|_{U^6_S}\|Q^S_2v^*_{N_2}\|_{U^6_S}\|w^*_{N_3}\|_{V^2_S}\notag \\
\lesssim& \delta^{\frac{1}{6}}N_1^{\frac{1}{2}}C(N_1,N_3)\|u_{N_1}\|_{U^2_S}\|v_{N_2}\|_{U^2_S}\|w_{N_3}\|_{V^2_S}.\nonumber
\end{align}
{\bf Region 5.} $|\xi_3|\gg |\xi_1-\xi_2|$ and $|\xi_3|\gg |\eta_1|\geq |\eta_2|$

This can be dealt with in the same way as Region 4.\vspace{3mm}

Next, we turn to the estimate of the nonlinear term.
\begin{prop} \label {Bilinear1} Let $\frac{5}{7}<s<1$ and $0<\delta\leq1$. We have
\begin{align}
&\|\int^t_0e^{(t-t')S}\chi(\frac{t}{\delta})(\partial_x+\partial_y)I(uv)(t')dt'\|_{Y^1}\lesssim \delta^{\frac{1}{6}}\|Iu\|_{Y^1}\|Iv\|_{Y^1} \label{estimate35}.
 \end{align}
\end{prop}
{\bf Proof.} From symmetry and definition of $Y^1$ we need to consider the following two terms
\begin{align}
J_1=\sum_{N_3}N_3^2\|\sum_{N_1\gg N_2}\int^t_0e^{(t-t')S}\chi(\frac{t}{\delta})(\partial_x+\partial_y)P_{N_3}I(u_{N_1}v_{N_2})(t')dt'\|^2_{U^2_S}\nonumber
\end{align}
and
\begin{align}
J_2=\sum_{N_3}N_3^2\|\sum_{N_1\sim N_2}\int^t_0e^{(t-t')S}\chi(\frac{t}{\delta})(\partial_x+\partial_y)P_{N_3}I(u_{N_1}v_{N_2})(t')dt'\|^2_{U^2_S}.\nonumber
\end{align}
Combining Proposition \ref{UV prop1} \rm(iii) with \eqref{estimate33}, we obtain
\begin{align}
J_1&\lesssim \sum_{N_3}N_3^2(\sum_{N_1\gg N_2}\mathop{sup}\limits_{\|w_{N_3}\|_{V^2_S}\leq1}|\int_{\Real^3}\chi(\frac{t}{\delta})I(u_{N_1}v_{N_2})
(\partial_x+\partial_y)w_{N_3}dxdydt|)^2   \notag \\
&\lesssim \delta\sum_{N_3}N_3^2(\sum_{N_1\gg N_2}N_2^{\frac{3}{2}-s}\|Iu_{N_1}\|_{U^2_S}\|Iv_{N_2}\|_{U^2_S})^2    \notag \\
&\lesssim \delta\sum_{N_1}N_1^2\|Iu_{N_1}\|^2_{U^2_S}\sum_{N_2}N_2^2\|Iv_{N_2}\|^2_{U^2_S}\sum_{N_2}N_2^{1-2s}\notag \\
&\lesssim \delta\|Iu\|^2_{Y^1}\|Iv\|^2_{Y^1}.\nonumber
 \end{align}
The second term can be controlled via \eqref{estimate34} and Minkowshi inequality,
\begin{align}
J_2&\lesssim\sum_{N_3}N_3^2(\sum_{N_1\sim N_2}\mathop{sup}\limits_{\|w_{N_3}\|_{V^2_S}\leq1}|\int_{\Real^3}\chi(\frac{t}{\delta})I(u_{N_1}v_{N_2})
(\partial_x+\partial_y)w_{N_3}dxdydt|)^2   \notag \\
   &\lesssim \delta^{\frac{1}{3}}(\sum_{N_3\sim N_1}N_3^2(N_1^{\frac{3}{2}-s}\|Iu_{N_1}\|_{U^2_S}\|Iv_{N_1}\|_{U^2_S})^2\notag \\
    & \ \ \ +\sum_{N_3\ll N}N_3^2(\sum_{N \lesssim N_1}N_1^{\frac{1}{2}}(\frac{N_1}{N})^{2-2s}\|Iu_{N_1}\|_{U^2_S} \|Iv_{N_1}\|_{U^2_S})^2 \notag \\
&\ \ \ +\sum_{N \lesssim N_3}N_3^2(\sum_{N_3\ll N_1}N_1^{\frac{1}{2}}N^{s-1}N_1^{2-2s}N_3^{s-1}\|Iu_{N_1}\|_{U^2_S}\|Iv_{N_1}\|_{U^2_S})^2) \notag \\
   &\lesssim \delta^{\frac{1}{3}}(\sum_{N_1}N_1^{5-2s}\|Iu_{N_1}\|_{U^2_S}^2\|Iv_{N_1}\|_{U^2_S}^2\notag \\
&\ \ \  +N^{4s-2}(\sum_{N \lesssim N_1}N_1^2N_1^{\frac{1}{2}-2s}\|Iu_{N_1}\|_{U^2_S}\|Iv_{N_1}\|_{U^2_S})^2 \notag \\
&\ \ \ +N^{2s-2}\sum_{N_3\ll N_1}N_3^{2s}(\sum_{N_1}N_1^{\frac{5}{2}-2s}\|Iu_{N_1}\|_{U^2_S}\|Iv_{N_1}\|_{U^2_S})^2) \notag \\
&\lesssim \delta^{\frac{1}{3}}(\sum_{N_1}N_1^4\|Iu_{N_1}\|_{U^2_S}^2\|Iv_{N_2}\|_{U^2_S}^2+(\sum_{N_1}N_1^2\|Iu_{N_1}\|_{U^2_S}\|Iv_{N_1}\|_{U^2_S})^2\notag \\
&\ \ \ +(\sum_{N_1}N_1^{\frac{5}{2}-s}\|Iu_{N_1}\|_{U^2_S}\|Iv_{N_1}\|_{U^2_S})^2) \notag \\
&\lesssim \delta^{\frac{1}{3}}\|Iu\|^2_{Y^1}\|Iv\|^2_{Y^1}.\nonumber
\end{align}

Hence, we complete the proof.

\begin{prop} \label{p2}
Let $\frac{5}{7}<s<1$. Assume $u_0$ satisfies $|E(Iu_0)|\leq 1$. Then there is a constant $\delta=\delta(\|u_0\|_{L^2(\Real^2)})$ and a unique solution $u$ to \eqref {ZK} on $[0,\delta]$,
 such that
$$\|Iu\|_{Y^1}\lesssim 1$$
 where the implicit constant is independent of $\delta$.
\end{prop}
{\bf Proof.} Acting multiplier operator $I$ on both sides of \eqref {ZK}, one can obtain
\begin{align}\partial_{t}Iu + (\partial^3_{x}+\partial^3_{y}) Iu +(\partial_{x}+\partial_{y}) I(u^2) = 0.\label{equationI}\end{align}
From Duhamel's principle, we get
$$Iu=\chi(\frac{t}{\delta})e^{tS}Iu_0-\int^t_0e^{(t-t')S}\chi(\frac{t}{\delta})(\partial_{x}+\partial_{y})I(u^2)(t')dt'.$$
Applying Proposition \ref {Bilinear1}, one obtain
\begin{align}
\|Iu\|_{Y^1} &\leq\|\chi(\frac{t}{\delta})e^{tS}Iu_0\|_{Y^1}+\|\int^t_0e^{(t-t')S}\chi(\frac{t}{\delta})(\partial_x+\partial_y)I(u^2)(t')dt'\|_{Y^1} \notag \\
             & \leq \|Iu_0\|_{H^1(\Real^2)} +C\delta^{\frac{1}{6}}\|Iu\|^2_{Y^1}, \nonumber
 \end{align}
and
$$\|Iu-Iv\|_{Y^1}\leq C\delta^{\frac{1}{6}}(\|Iu\|_{Y^1}+\|Iv\|_{Y^1})\|Iu-Iv\|_{Y^1}.$$
Then, the contraction mapping principle ensure the existence of local solution. Moreover, a bootstrap argument yields $\|Iu\|_{Y^1}\lesssim 1$.

According to equation \eqref {equationI} and integration by parts, we have
\begin{align}
\frac{dE(Iu)(t)}{dt}& =- \int_{\Real^2} (\Delta Iu -\partial_{x}\partial_{y}Iu+(Iu)^2)\partial_{t}Iudxdy   \notag \\
& = \int_{\Real^2}((\partial^3_{x}+\partial^3_{y})Iu+(\partial_{x}+\partial_{y})Iu^2)((\Delta-\partial_{x}\partial_{y})Iu +(Iu)^2)dxdy   \notag \\
& = \int_{\Real^2}(\partial^3_{x}+\partial^3_{y})Iu((Iu)^2-Iu^2)dxdy+\int_{\Real^2}(\partial_{x}+\partial_{y})Iu^2(Iu)^2dxdy.   \nonumber
\end{align}
Integrating in time and applying the Parseval formula provide
\begin{align}
E(Iu)(\delta)- E(Iu)(0)& =\int^{\delta}_0\int_{\sum^3_{j=1}\zeta_j=0}(\xi^3_1+\eta^3_1)(1-\frac{m(\zeta_2+\zeta_3)}{m(\zeta_2)m(\zeta_3)})\prod^3_{j=1}\widehat{Iu}(\zeta_j)
 \notag \\
& +\int^{\delta}_0\int_{\sum^4_{j=1}\zeta_j=0}(\xi_1+\xi_2+\eta_1+\eta_2)\frac{m(\zeta_1+\zeta_2)}{m(\zeta_1)m(\zeta_2)}\prod^4_{j=1}\widehat{Iu}(\zeta_j). \label{In1}
\end{align}

The following preliminaries are effective for controlling the growth of $E(Iu)(t)$.
\begin{lem} \label {Enery lema}  Let $N_1, N_2, N_3$ be dyadic numbers. Then, for all $0<\delta\leq1$ we have
\\if $N_1\sim N_2\sim N_3\gtrsim N$,
\begin{align}
&|\int_{\sum^3_{j=1}\lambda_j=0}(\xi^3_1+\eta^3_1)(1-\frac{m(\zeta_2+\zeta_3)}{m(\zeta_2)m(\zeta_3)})\prod^3_{j=1}\F (\chi_{\delta}u_j)(\lambda_j)|\notag \\
\lesssim& \delta^{\frac{1}{6}}N^{s-1}N_1^{\frac{7}{2}-s}\prod^3_{j=1}\|u_j\|_{U^2_S} \label{estimateE1},\end{align}
if $N_1\sim  N_2\gg N_3$, $N_1\gtrsim N$,
\begin{align}
&|\int_{\sum^3_{j=1}\lambda_j=0}(\xi_1\xi_2\xi_3+\eta_1\eta_2\eta_3)\prod^3_{j=1}\F (\chi_{\delta}u_j)(\lambda_j)|\notag \\
\lesssim& \delta^{\frac{1}{2}}N_1N_3^{\frac{3}{2}}\prod^3_{j=1}\|u_j\|_{U^2_S} \label{estimateE2},
\end{align}
and
\begin{align}
&|\int_{\sum^3_{j=1}\lambda_j=0}\frac{\sum^3_{j=1}(\xi^3_j+\eta^3_j)m^2(\zeta_j)}{m(\zeta_1)m(\zeta_2)m(\zeta_3)}\prod^3_{j=1}\F (\chi_{\delta}u_j)(\lambda_j)|\notag \\
\lesssim&\delta^{\frac{1}{2}}N_1N_3^{\frac{3}{2}}(1\lor(\frac{N_3}{N})^{1-s})\prod^3_{j=1}\|u_j\|_{U^2_S}. \label{estimateE3}
\end{align}
\end{lem}
{\bf Proof.} When $N_1\sim N_2\sim N_3\gtrsim N$,
$$|1-\frac{m(\zeta_2+\zeta_3)}{m(\zeta_2)m(\zeta_3)}|\lesssim\frac{m(\zeta_1)}{m(\zeta_2)m(\zeta_3)}\lesssim (\frac{N_1}{N})^{1-s}$$
and 
$$|\xi^3_1+\eta^3_1|\lesssim |\xi_1+\eta_1|N^2_1.$$ 
Then one can get \eqref{estimateE1} as the second part of Lemma \ref {local lemma2}.

In order to prove \eqref{estimateE2}, we assume $|\xi_1\xi_2\xi_3|\geq|\eta_1\eta_2\eta_3|$ by symmetry. \\
{\bf Case 1.} $|\xi_2|\lesssim|\xi_1-\xi_3|$

From \eqref{estimate24}, we have 
\begin{align}
&|\int_{\sum^3_{j=1}\lambda_j=0}(\xi_1\xi_2\xi_3+\eta_1\eta_2\eta_3)\prod^3_{j=1}\F (\chi_{\delta}u_j)(\lambda_j)|\notag \\
\lesssim&N_1N_3|\int_{\sum^3_{j=1}\lambda_j=0}|\xi_2|^{\frac{1}{2}}|\xi_1-\xi_3|^{\frac{1}{2}}\prod^3_{j=1}\F (\chi_{\delta}u_j)(\lambda_j)|\notag \\
\lesssim &N_1N_3\|\int|\xi_1+\xi_3|^{\frac{1}{2}}|\xi_1-\xi_3|^{\frac{1}{2}}\tilde u_{1}(\lambda_1)\tilde u_{3}(\lambda_3)\|_{L^2(\Real^3)}\|\chi(\frac{t}{\delta})u_{2}\|_{L^2(\Real^3)}\notag \\
\lesssim& \delta^{\frac{1}{2}}N_1N_3^{\frac{3}{2}}\prod^3_{j=1}\|u_j\|_{U^2_S}. \nonumber
\end{align}
{\bf Case 2.} $|\xi_2|\gg|\xi_1-\xi_3|$

In this case, $|\xi_2|\sim|\xi_1|\sim|\xi_3|\lesssim N_3\ll N_1\sim N_2$. Hence, $$|\xi_1\xi_2\xi_3|\lesssim N_3^{\frac{5}{2}}|\xi_1|^{\frac{1}{8}}|\eta_1|^{\frac{1}{8}}|\xi_2|^{\frac{1}{8}}|\eta_2|^{\frac{1}{8}}.$$
From \eqref{estimate23}, we have 
\begin{align}
&|\int_{\sum^3_{j=1}\lambda_j=0}(\xi_1\xi_2\xi_3+\eta_1\eta_2\eta_3)\prod^3_{j=1}\F (\chi_{\delta}u_j)(\lambda_j)|\notag \\
\lesssim&|\int_{\sum^3_{j=1}\lambda_j=0}|\xi_1\xi_2\xi_3|\prod^3_{j=1}\F (\chi_{\delta}u_j)(\lambda_j)|\notag \\
\lesssim&N_3^{\frac{5}{2}}|\int_{\sum^3_{j=1}\lambda_j=0}|\xi_1|^{\frac{1}{8}}|\eta_1|^{\frac{1}{8}}|\xi_2|^{\frac{1}{8}}|\eta_2|^{\frac{1}{8}}\prod^3_{j=1}\F (\chi_{\delta}u_j)(\lambda_j)|\notag \\
\lesssim &N_3^{\frac{5}{2}}\|I^{\frac{1}{8}}_xI^{\frac{1}{8}}_yu_{1}\|_{L^4(\Real^3)}\|I^{\frac{1}{8}}_xI^{\frac{1}{8}}_yu_{2}\|_{L^4(\Real^3)}\|\chi(\frac{t}{\delta})u_{3}\|_{L^2(\Real^3)}\notag \\
\lesssim& \delta^{\frac{1}{2}}N_1N_3^{\frac{3}{2}}\prod^3_{j=1}\|u_j\|_{U^2_S}. \nonumber
\end{align}

It's easy to know that $N_1\sim  N_2\gg N_3$ implies $|\xi_1|\sim |\xi_2|\sim  N_1$ or $|\eta_1|\sim |\eta_2|\sim  N_1$. Without loss of generality one can $|\xi_1|\sim |\xi_2|\sim  N_1$, then mean-value theorem gives the bound of the symbol
$$\sum^3_{j=1}(\xi^3_j+\eta^3_j)m^2(\zeta_j)\lesssim N_1^2 \mathop{max}\left\{|\xi_3|, |\eta_3|\right\}m^2(\zeta_1).$$
Hence, one can show \eqref{estimateE3} in a similar way as above. 

We complete the proof of this lemma.

\begin{prop} \label{p3}
Let $\frac{5}{7}<s<1$. We have
\begin{align}
 |\int^{\delta}_0 \Lambda_3( (\xi^3_1+\eta^3_1)(1-\frac{m(\zeta_2+\zeta_3)}{m(\zeta_2)m(\zeta_3)});Iu)|\lesssim
  \delta^{\frac{1}{6}}N^{-\frac{1}{2}+}\|Iu\|^3_{Y^1} .\label{Lambda(3)}
\end{align}
\end{prop}
    {\bf Proof.}  The fact $$\sum_{N_1}N_1^{0-}N_1\|P_{N_1}Iu\|_{U^2_S}\leq (\sum_{N_1}N_1^{0-})^{\frac{1}{2}}\|Iu\|_{Y^1}$$ indicates that it's suffices to prove    
\begin{align}
| & \int^{\delta}_0\int_{\sum^3_{j=1}\zeta_j=0}[(\xi^3_1+\eta^3_1)(1-\frac{m(\zeta_2+\zeta_3)}{m(\zeta_2)m(\zeta_3)})]_{sym}\prod^3_{j=1}\widehat{u_j}(\zeta_j)| \notag \\
\lesssim &\delta^{\frac{1}{6}}N^{-\frac{1}{2}+}(N_1N_2N_3)^{1-}\prod^3_{j=1}\|u_j\|_{U^2_S} \label{Lambda(3.1)}
\end{align}
 for any function $u_j$ ($j=1,2,3$) with frequencies supported on $|\zeta_j|\sim  N_j$.
 
We denote $L_1$ the left hand of \eqref{Lambda(3.1)} and $M(\zeta_1,\zeta_2,\zeta_3)=[(\xi^3_1+\eta^3_1)(1-\frac{m(\zeta_2+\zeta_3)}{m(\zeta_2)m(\zeta_3)})]_{sym}$. Furthermore, one can assume $N_1\sim N_2\gtrsim N_3$ and $N_1\gtrsim N$.\\
{\bf Case 1.} $N_1\sim N_2\sim N_3\gtrsim N$

We obtain by \eqref{estimateE1}
\begin{align}
 L_1 &\lesssim \delta^{\frac{1}{6}}N^{s-1}N_1^{\frac{7}{2}-s}\prod^3_{j=1}\|u_j\|_{U^2_S}  \notag \\
&\lesssim \delta^{\frac{1}{6}}N^{s-1}N_1^{\frac{1}{2}-s+}N_1^{3-}\prod^3_{j=1}\|u_j\|_{U^2_S}\notag \\
&\lesssim \delta^{\frac{1}{6}}N^{-\frac{1}{2}+} N_1^{3-}\prod^3_{j=1}\|u_j\|_{U^2_S}. \nonumber
\end{align}
{\bf Case 2.} $N_1\sim N_2\gg N_3, N_1 \gtrsim N $

We write $M=M_1-M_2$, where $M_1=\sum^3_{j=1}(\xi^3_j+\eta^3_j)$ and
$M_2=\frac{\sum^3_{j=1}(\xi^3_j+\eta^3_j)m^2(\zeta_j)}{m(\zeta_1)m(\zeta_2)m(\zeta_3)}$.
The corresponding terms are
\begin{align}L_{1,1}=|\int^{\delta}_0\int_{\sum^3_{j=1}\zeta_j=0}M_1(\zeta_1,\zeta_2,\zeta_3)\prod^3_{j=1}\widehat{u_j}\zeta_j)| \nonumber \end{align}
and \begin{align}L_{1,2}=|\int^{\delta}_0\int_{\sum^3_{j=1}\zeta_j=0}M_2(\zeta_1,\zeta_2,\zeta_3)\prod^3_{j=1}\widehat{u_j}(\zeta_j)|. \nonumber \end{align}
Estimate \eqref{estimateE2} tells us that
\begin{align}
L_{1,1}&\lesssim \delta^{\frac{1}{2}}N_1N_3^{\frac{3}{2}}N_1^{-2+}N_3^{-1+}\prod^3_{j=1}N_j^{1-}\|u_j\|_{U^2_S}
\notag \\
&\lesssim \delta^{\frac{1}{2}}N^{-\frac{1}{2}+}\prod^3_{j=1}N_j^{1-}\|u_j\|_{U^2_S}.
\nonumber
 \end{align}
Using \eqref{estimateE3}, we get the contribution of $L_{1,2}$
\begin{align}L_{1,2}&\lesssim \delta^{\frac{1}{2}}N_1N_3^{\frac{3}{2}}(1\lor(\frac{N_3}{N})^{1-s})N_1^{-2+}N_3^{-1+}\prod^3_{j=1}N_j^{1-}\|u_j\|_{U^2_S} \notag \\
&\lesssim \delta^{\frac{1}{2}}N^{-\frac{1}{2}+}\prod^3_{j=1}N_j^{1-}\|u_j\|_{U^2_S}.
\nonumber
\end{align}
This complete the proof of \eqref{Lambda(3.1)}, and hence \eqref{Lambda(3)}.

\begin{prop} \label{p4}
Let $\frac{5}{7}<s<1$. We have
\begin{align}
 |\int^{\delta}_0 \Lambda_4((\xi_1+\xi_2+\eta_1+\eta_2)\frac{m(\zeta_1+\zeta_2)}{m(\zeta_1)m(\zeta_2)};Iu)|\lesssim \delta^{\frac{1}{6}}N^{-1+}\|Iu\|^4_{Y^1}\label{Lambda(4)}.
\end{align}
\end{prop}
{\bf Proof.} As previous discussion, it suffices to show
\begin{align}
| & \int^{\delta}_0 \int_{\sum^4_{j=1}\lambda_j=0}[(\xi_1+\xi_2+\eta_1+\eta_2)\frac{m(\zeta_1+\zeta_2)}{m(\zeta_1)m(\zeta_2)}]_{sym}\prod^4_{j=1}\F(\chi_{\delta}u_j)(\lambda_j)| \notag \\
\lesssim &\delta^{\frac{1}{6}}N^{-1+}\prod^4_{j=1}N_j^{1-}\|u_j\|_{U^2_S}. \label{Lambda(4.1)}
\end{align}

If $N \gg N_j$ $(j=1,2,3,4)$, then it is easy to see $[(\xi_1+\xi_2)\frac{m(\zeta_1+\zeta_2)}{m(\zeta_1)m(\zeta_2)}]_{sym}=0$, \eqref {Lambda(4.1)} holds trivially. We may assume $N_1\geq N_2, N_3\geq N_4$ by  symmetry and let  $L_2$ denote the left hand of \eqref{Lambda(4.1)}.

From H\"{o}rder's inequalities and \eqref{estimate22}, we obtain
\begin{align}
L_2&\lesssim \|\int^{\delta}_0 \int_{\zeta=\zeta_1+\zeta_2}(\xi+\eta)\frac{m(\zeta)}{m(\zeta_1)m(\zeta_2)}\widehat u_1(\zeta_1)\widehat u_2(\zeta_2)d\zeta_1dt\|_{L^2}\|\chi_{\delta}u_3\|_{L^4}\|\chi_{\delta}u_4\|_{L^4}\notag \\
&\lesssim \delta^{\frac{1}{3}}N_1\|\F^{-1}\frac{\widehat u_1}{m(\zeta_1)}\|_{L^6_tL^4_{x,y}}\|\F^{-1}\frac{\widehat u_2}{m(\zeta_2)}\|_{L^6_tL^4_{x,y}}\|u_3\|_{L^6_tL^4_{x,y}}\| u_4\|_{L^6_tL^4_{x,y}}\notag \\
&\lesssim \delta^{\frac{1}{3}}N_1\|\F^{-1}\frac{\widehat u_1}{m(\zeta_1)}\|_{U^6_S}\|\F^{-1}\frac{\widehat u_2}{m(\zeta_2)}\|_{U^6_S}\|u_3\|_{U^6_S}\| u_4\|_{U^6_S}\notag \\
&\lesssim \delta^{\frac{1}{3}}N_1\|\F^{-1}\frac{\widehat u_1}{m(\zeta_1)}\|_{U^2_S}\|\F^{-1}\frac{\widehat u_2}{m(\zeta_2)}\|_{U^2_S}\|u_3\|_{U^2_S}\| u_4\|_{U^2_S}\notag \\
&\lesssim \delta^{\frac{1}{3}}\frac{N_1}{m(N_1)m(N_2)}\prod^4_{j=1}\|u_j\|_{U^2_S}. \nonumber
\end{align}
{\bf Case 1.} $N_1\ll N_3$

This implies $N_3\sim N_4\gg N_1\geq N_2$,  $N_3\sim N_4\gtrsim N.$ Combining the  pointwise bound about the symbol

$$ \frac{N^{-1+}_j}{m(N_j)}\lesssim 1$$
for all $N_j$, we have
\begin{align}
L_2&\lesssim \delta^{\frac{1}{3}}N_1N^{-2+}_3\prod^4_{j=1}N_j^{1-}\|u_j\|_{U^2_S}\notag \\
&\lesssim \delta^{\frac{1}{3}}N^{-1+}\prod^4_{j=1}N_j^{1-}\|u_j\|_{U^2_S}. \nonumber
\end{align}
{\bf Case 2.} $N_1\gtrsim N_3, N_1 \sim N_2\gtrsim N$

Similarly to Case 1, one obtain that

\begin{align}
  L_2&\lesssim \delta^{\frac{1}{3}}\frac{N_1}{m(N_1)m(N_2)}N_1^{-2+}\prod^4_{j=1}N_j^{1-}\|u_j\|_{U^2_S}\notag \\
  &\lesssim \delta^{\frac{1}{3}}N_1(\frac{N_1}{N})^{2(1-s)}N_1^{-2+}\prod^4_{j=1}N_j^{1-}\|u_j\|_{U^2_S}\notag \\
&\lesssim \delta^{\frac{1}{3}}N^{-1+}\prod^4_{j=1}N_j^{1-}\|u_j\|_{U^2_S}. \nonumber
\end{align}
{\bf Case 3.} $N_1\gtrsim N_3, N_1 \gg N_2$

In this case, $N_1 \sim N_3\gtrsim N$, we control $L_2$ by using \eqref {estimate31},
\begin{align}
  T_2&\lesssim  \|\int_{\lambda=\lambda_1+\lambda_2}(\xi+\eta)\frac{m(\zeta_1+\zeta_2)}{m(\zeta_1)m(\zeta_2)}\tilde u_1(\lambda_1)\tilde u_2(\lambda_2)d\lambda_1\|_{L^2_{\lambda}}\|\chi_{\delta}u_3\|_{L^4}\|\chi_{\delta}u_4\|_{L^4}\notag \\
  &\lesssim \delta^{\frac{1}{6}}\frac{N_2^{\frac{1}{2}}}{m(N_2)}\prod^4_{j=1}\|u_j\|_{U^2_S}\notag \\
    &\lesssim \delta^{\frac{1}{6}}N_2^{\frac{1}{2}}(1\lor (\frac{N_2}{N})^{1-s})N_1^{-2+}N_2^{-1+}\prod^4_{j=1}N_j^{1-}\|u_j\|_{U^2_S}\notag \\
     &\lesssim \delta^{\frac{1}{6}}N^{-2+}\prod^4_{j=1}N_j^{1-}\|u_j\|_{U^2_S}. \nonumber
\end{align} 
Therefore we complete the proof of \eqref{Lambda(4.1)}.

It's not hard to obtain the following proposition in a similar way.
\begin{prop} \label{p5}
Let $\frac{5}{7}<s<1$. We have
\begin{align}
 |\int^{\delta}_0 \Lambda_3((\xi_3+\eta_3)\frac{m(\zeta_1+\zeta_2)}{m(\zeta_1)m(\zeta_2)};Iu)|\lesssim \delta^{\frac{2}{3}}N^{-2+}\|Iu\|^3_{Y^1}\label{Lambda(32)}.
\end{align}
\end{prop}
\begin{prop} \label {Prop E incresement}
Let $\frac{5}{7}<s<1, N\gg 1$. Assume $u_0$ satisfies $|E(Iu_0)|\leq 1$. Then there is a constant $\delta=\delta(\|u_0\|_{L^2(\Real^2)})>0$ so that there exists a unique solution
$$u(x,y,t)\in C([0,\delta],H^s(\Real^2))$$
 of \eqref {ZK} satisfying
\begin{align}
E(I_Nu)(\delta)=E(I_Nu)(0)+\delta^{\frac{1}{6}}O(N^{-\frac{1}{2}+}). \label {E incresement}
 \end{align}
\end{prop}
{\bf Proof.} From Proposition \ref{p2}, there exsits a unique solution $u$ to \eqref {ZK} on $[0,\delta]$ satisfying $\|Iu\|_{Y^1}\lesssim 1$.

Combining \eqref {In1}, Proposition \ref{p3} and Proposition \ref {p4}, one has
\begin{align}
 |E(Iu)(\delta)- E(Iu)(0)|& =|\int^{\delta}_0 \Lambda_3( (\xi^3_1+\eta^3_1)(1-\frac{m(\zeta_2+\zeta_3)}{m(\zeta_2)m(\zeta_3)});Iu)\notag \\
 &+\int^{\delta}_0 \Lambda_4((\xi_1+\xi_2+\eta_1+\eta_2)\frac{m(\zeta_1+\zeta_2)}{m(\zeta_1)m(\zeta_2)};Iu)|\notag \\
 & \lesssim \delta^{\frac{1}{6}}N^{-\frac{1}{2}+}\|Iu\|^3_{Y^1}+\delta^{\frac{1}{6}}N^{-1+}\|Iu\|^4_{Y^1}\notag \\
 & \lesssim \delta^{\frac{1}{6}}N^{-\frac{1}{2}+}.\nonumber
\end{align}

Our purpose is to construct a solution $u$ on $[0,T]$ for any $T>0$. Note that $u$ is a solution to \eqref {ZK}, then $u_\lambda(x,y,t)=\lambda^2u(\lambda x,\lambda y,\lambda^3t)$ is also the solution to \eqref {ZK}. Hence it suffices to acquire the well-posedness for $u_\lambda$ on $[0,\frac{T}{\lambda^3}]$.

From Proposition 3.1 the energy will be arbitrarily small by taking $\lambda$ small,
\begin{align}
E(I_Nu_{\lambda,0})&\lesssim N^{2(1-s)}\|u_{\lambda,0}\|^2_{\dot H^s(\Real^2)}+\|u_{\lambda,0}\|^3_{L^3(\Real^2)}\notag \\
 & \lesssim N^{2(1-s)}\lambda^{2(s+1)}\|u_0\|^2_{\dot H^s(\Real^2)}+\lambda^4\|u_0\|^3_{L^3(\Real^2)}\notag \\
 & \lesssim (N^{2(1-s)}\lambda^{2(s+1)}+\lambda^4)(1+\|u_0\|_{H^s(\Real^2)})^3.\nonumber
\end{align}
Assume $N\gg1$ is given ($N$ will be chose later), then we have 
$$E(I_Nu_{\lambda,0})\leq \frac{1}{4}$$
by setting
\begin{align}
\lambda=\lambda(N,\|u_0\|_{H^s(\Real^2)})\sim N^{\frac{s-1}{s+1}}.\nonumber
\end{align}
Now applying Proposition \ref {Prop E incresement} to $u_{\lambda,0}$, one gets
$$E(Iu_\lambda)(\delta)\leq\frac{1}{4}+CN^{-\frac{1}{2}+}<\frac{1}{2}.$$
Thus from Proposition \ref {p2} the solution $u_\lambda$ can be extended to $t\in[0,2\delta]$. Iterating this procedure $M$ steps, we have 
$$E(Iu_\lambda)(t)\leq\frac{1}{4}+CMN^{-\frac{1}{2}+}$$ 
for $t\in[0,(M+1)\delta]$. That's to say, as long as $MN^{-\frac{1}{2}+}\lesssim1$ the solution $u_\lambda$ can be extended to $t\in[0, (M+1)\delta].$
Taking $N(T)\sim T^{\frac{2(s+1)}{7s-5}+}\gg1$, then
$$(M+1)\delta \sim N^{\frac{1}{2}-}\delta \sim TN^{\frac{3(1-s)}{s+1}}\sim \frac{T}{\lambda^3}.$$
Hence, we obtain global well-posedness for \eqref {ZK} when  $s>\frac{5}{7}$. Moreover, from Proposition 3.1 we have
\begin{align}\|u(T)\|_{ H^s(\Real^2)}&\lesssim \lambda^{-s-1}\|u_\lambda(\frac{T}{\lambda^3})\|_{ H^s(\Real^2)}\notag \\
&\lesssim \lambda^{-s-1}(|E(Iu_\lambda)(\frac{T}{\lambda^3})|^{\frac{1}{2}}+\|u_{\lambda,0}\|_{L^2(\Real^2)}+\|u_{\lambda,0}\|^2_{L^2(\Real^2)})\notag \\
&\lesssim \lambda^{-s-1} (1+\|u_0\|_{ H^s(\Real^2)})^2\notag \\
&\lesssim T^{\frac{2(1-s)(1+s)}{7s-5}-} (1+\|u_0\|_{ H^s(\Real^2)})^2.\nonumber\end{align}

\section{Global attractor}

This section follows Tsugawa's idea which he applied to show the existence of global attarctor for KdV equation on Sobolev spaces of negative index (see \cite{Tsugawa} ).

We find a representation of the rescaled equation associated to \eqref{dfZK}
\begin{equation}
  \left\{
   \begin{aligned}
   &\partial_{t}v +( \partial^3_{x}+\partial^3_{y}) v +( \partial_{x}+\partial_{y}) v^2+\gamma \lambda^{-3}v =\lambda^{-3}g,     \quad \\
   &v(x,y,0)=v_0(x,y)\in H^s(\Real^2), \label{rdfZK} \\
   \end{aligned}
   \right.
\end{equation}
with $v(x,y,t)=\lambda^{-2}u(\lambda^{-1}x,\lambda^{-1}y,\lambda^{-3}t)$, $v_0(x,y)=\lambda^{-2}u_0(\lambda^{-1}x,\lambda^{-1}y)$ and $g(x,y)=\lambda^{-2}f(\lambda^{-1}x,\lambda^{-1}y)$.

From the definition of rescaled operator $I'$, one clearly knows
\begin{align}
 \ \ \ \|I'v\|_{L^2}= \lambda^{-1}\|Iu\|_{L^2} \ , \ \ \ \|I'g\|_{L^2}= \lambda^{-1}\|If\|_{L^2},
\nonumber
\end{align}
and
\begin{align}
 \ \ \ \|I'v\|_{\dot H^1}= \lambda^{-2}\|Iu\|_{\dot H^1} \ , \ \ \ \|I'g\|_{\dot H^1}= \lambda^{-2}\|If\|_{\dot H^1}. 
\nonumber
\end{align}

Now, we give the time local result for weakly damped forced ZK equation.
\begin{prop} \label {Prop lwp dfZK}
Let $\frac{5}{7}<s<1$. Assume  $I'v_0 \in H^1(\Real^2) $ and $I'g \in H^1(\Real^2) $, then there is a constant $\delta=\delta(\|I'v_0 \|_{H^1(\Real^2)}, \lambda^{-3}\|I'g \|_{H^1(\Real^2)},\gamma\lambda^{-3})>0$ so that there exists a unique solution
$$v(x,y,t)\in C([0,\delta],H^s(\Real^2))$$
 of \eqref {rdfZK} satisfying
\begin{align}
  \|I'v \|_{Y^1}\lesssim \|I'v_0 \|_{H^1(\Real^2)}+ \lambda^{-3}\|I'g \|_{H^1(\Real^2)},\nonumber
   \end{align}
  and 
  \begin{align}
  \mathop{sup}\limits_{t\in[0,\delta]} \|I'v(t)\|_{H^1}\lesssim \|I'v_0 \|_{H^1(\Real^2)}+ \lambda^{-3}\|I'g \|_{H^1(\Real^2)}.\nonumber
 \end{align}
\end{prop}
{\bf Proof.} Acting $I'$ on \eqref {rdfZK} gives
\begin{align}
  \partial_{t}I'v +( \partial^3_{x}+\partial^3_{y})I'v +( \partial_{x}+\partial_{y})I'v^2+\gamma \lambda^{-3}I'v =\lambda^{-3}I'g \label{I'rdfZK}
\end{align}
which also can been written as an intergral equation
$$I'v=\mathscr{T}I'v$$
where
$$\mathscr{T}I'v=\chi(\frac{t}{\delta})e^{tS}Iv_0-\int^t_0e^{(t-t')S}\chi(\frac{t}{\delta})((\partial_{x}+\partial_{y})I'v^2+\gamma \lambda^{-3}I'v -\lambda^{-3}I'g)dt'.$$
By the duality of $U^p$, we have
\begin{align}
  &\|\int^t_0e^{(t-t')S}\chi(\frac{t}{\delta})(\gamma \lambda^{-3}I'v -\lambda^{-3}I'g)dt' \|_{Y^1}\notag \\
  \lesssim& (\sum_{N_1}N_1^2(\mathop{sup}\limits_{\|w\|_{V^2_S}\leq1}|\int_{\Real^3}\chi(\frac{t}{\delta})(\gamma \lambda^{-3}P_{N_1}I'v -\lambda^{-3}P_{N_1}I'g)wdxdydt|)^2)^{\frac{1}{2}} \notag \\
  \lesssim& (\sum_{N_1}N_1^2(\mathop{sup}\limits_{\|w\|_{V^2_S}\leq1}\|\gamma \lambda^{-3}P_{N_1}I'v -\lambda^{-3}P_{N_1}I'g\|_{L^{\infty}_tL^2_{x,y}}\|\chi(\frac{t}{\delta}) w\|_{L^1_tL^2_{x,y}})^2)^{\frac{1}{2}} \notag \\
 \lesssim& \delta(\sum_{N_1}N_1^2(\|\gamma \lambda^{-3}P_{N_1}I'v\|_{U^2_S}+\lambda^{-3}\|P_{N_1}I'g\|_{L^2_{x,y}})^2)^{\frac{1}{2}} \notag \\
 \lesssim& \delta(\gamma \lambda^{-3}\|I'v\|_{Y^1}+\lambda^{-3}\|I'g\|_{H^1}). \label{lwp dfterm} 
\end{align}
Set $$B=\{I'v\in Y^1 \big|\  \|I'v\|_{Y^1}<C_0(\|I'v_0\|_{H^1}+\lambda^{-3}\|I'g\|_{H^1})\}.$$
Applying Proposition \ref {Bilinear1} and \eqref{lwp dfterm}, on $B$ one gets
\begin{align}
  \|\mathscr{T}I'v\|_{Y^1}&\lesssim \|I'v_0 \|_{H^1}+\delta^{\frac{1}{6}}\|I'v\|^2_{Y^1}+ \delta\gamma \lambda^{-3}\|I'v\|_{Y^1}+\delta\lambda^{-3}\|I'g\|_{H^1}\notag \\
&\lesssim \delta^{\frac{1}{6}}C^2_0(\|I'v_0\|_{H^1}+\lambda^{-3}\|I'g\|_{H^1})^2+(1+C_0\delta\gamma \lambda^{-3})(\|I'v_0\|_{H^1}+\lambda^{-3}\|I'g\|_{H^1}), \nonumber
\end{align}
and
\begin{align}
  \|\mathscr{T}I'v_1-\mathscr{T}I'v_2\|_{Y^1}&\lesssim (\delta\gamma \lambda^{-3}+\delta^{\frac{1}{6}}\|I'v_1\|_{Y^1}+\delta^{\frac{1}{6}}\|I'v_2\|_{Y^1})\|I'v_1-I'v_2\|_{Y^1}\notag \\
  &\lesssim (\delta\gamma \lambda^{-3}+\delta^{\frac{1}{6}}C_0(\|I'v_0\|_{H^1}+\lambda^{-3}\|I'g\|_{H^1}))\|I'v_1-I'v_2\|_{Y^1}.\nonumber
\end{align}
If we assume
\begin{align}
   \lambda^{-3}\gamma\ll1 \ , \ \ \ \ \|I'v_0\|_{H^1}\ll1 \ , \ \ \ \ \lambda^{-3}\|I'g\|_{H^1}\ll1 \ , \label{lwp condi} 
\end{align}
then
$$\mathscr{T}: B\rightarrow B$$
is a strict contraction mapping. 

Finally, fixing $\lambda=\lambda_0$ and $I'v_0$, $I'g\in H^1$, we consider $\sigma-$scaling of $v$
$$w(x,y,t)=\sigma^{-2}v(\sigma^{-1}x,\sigma^{-1}y,\sigma^{-3}t).$$
It's equivalent to consider well-posedness on $[0, \sigma^3\delta]$ for $w$.

Observe that
$$ (\sigma\lambda_0)^{-3}\gamma\ll1,$$
$$\|I''w_0\|_{H^1}\lesssim \sigma^{-1}\|I'v_0\|_{H^1}\ll1,$$
$$(\sigma\lambda_0)^{-3}\|I'g\|_{H^1}\ll1$$
provided $\sigma$ is chosen to be sufficiently large, which verifies \eqref{lwp condi}. It means that Cauchy problem for $w$ is well-posed on the time interval $[0,1]$. Hence, \eqref{I'rdfZK} is locally well-posed on $[0,\sigma^{-3}]$. We complete the proof.

In the next place, we explore the increment of $I'v$ through modified energy $E(I'v)$.

From \eqref{I'rdfZK}, we obtain
\begin{align}
\frac{dE(I'v)(t)}{dt}& =- \int_{\Real^2}(\Delta I'v -\partial_{x}\partial_{y}I'v+(I'v)^2)\partial_{t}I'vdxdy   \notag \\
& = \int_{\Real^2}((\partial^3_{x}+\partial^3_{y})I'v+(\partial_{x}+\partial_{y})I'v^2)((\Delta-\partial_{x}\partial_{y})I'v +(I'v)^2 )dxdy \notag \\
&\ \ \ +\int_{\Real^2}(\gamma \lambda^{-3}I'v-\lambda^{-3}I'g)((\Delta-\partial_{x}\partial_{y})I'v +(I'v)^2 )dxdy \notag \\
& = -2\gamma \lambda^{-3}E(I'v)-\int_{\Real^2}\lambda^{-3}I'g((\Delta-\partial_{x}\partial_{y})I'v +(I'v)^2 )dxdy \notag \\
&\ \ \ +\frac{1}{3}\gamma \lambda^{-3}\int_{\Real^2}(I'v)^3dxdy+\Lambda_3( (\xi^3_1+\eta^3_1)(1-\frac{m(\zeta_2+\zeta_3)}{m(\zeta_2)m(\zeta_3)});I'v)\notag \\
&\ \ \ +\Lambda_4((\xi_1+\xi_2+\eta_1+\eta_2)\frac{m(\zeta_1+\zeta_2)}{m(\zeta_1)m(\zeta_2)};I'v).   \nonumber
\end{align}
This implies
\begin{align}
  \frac{d}{dt}E(I'v)(t)e^{2\gamma \lambda^{-3}t}& =-\int_{\Real^2}\lambda^{-3}I'g((\Delta-\partial_{x}\partial_{y})I'v +(I'v)^2 )e^{2\gamma \lambda^{-3}t}dxdy \notag \\
&\ \ \ +\frac{1}{3}\gamma \lambda^{-3}\|I'v\|^3_{L^3}e^{2\gamma \lambda^{-3}t}+(\Lambda_3+\Lambda_4)e^{2\gamma \lambda^{-3}t}. \label{E(I'v)exp}
\end{align}
Integrating \eqref{E(I'v)exp} over $[0,T']$, one gets
\begin{align}
 & E(I'v)(T')e^{2\gamma \lambda^{-3}T'}-E(I'v)(0) \notag \\
= &-\int^{T'}_0\int_{\Real^2}\lambda^{-3}I'g[(\Delta-\partial_{x}\partial_{y})I'v+(I'v)^2 ]e^{2\gamma \lambda^{-3}t}dxdydt\notag \\
&+\frac{1}{3}\gamma \lambda^{-3}\int^{T'}_0\|I'v\|^3_{L^3}e^{2\gamma \lambda^{-3}t}dt+\int^{T'}_0(\Lambda_3+\Lambda_4)e^{2\gamma \lambda^{-3}t}dt.  \label{E(I'v)}
\end{align}

\begin{lem} \label{I'v lemma1}
  Assume that $v$ is a solution of  \eqref{rdfZK} on $[0,T']$.
Then, we have
\begin{align}
  &\mathop{sup}\limits_{t\in[0,T']}\|I'v(t)\|^2_{L^2}e^{2\gamma \lambda^{-3}t} \notag \\
\leq& C_1  (\|I'v_0\|^2_{L^2}+\frac{1}{\gamma^2}\|I'g\|^2_{L^2}e^{2\gamma \lambda^{-3}T'}+|\int^{T'}_0e^{2\gamma \lambda^{-3}t}\Lambda'_3dt|),  \label{lemma1 2}
\end{align}
and
\begin{align}
  &\mathop{sup}\limits_{t\in[0,T']}\|I'v(t)\|^2_{\dot H^1}e^{2\gamma \lambda^{-3}t} \notag \\
  \leq& C_1(\|I'v_0\|^2_{\dot H^1}+\|I'v_0\|^4_{L^2}+\frac{1}{\gamma^2} \|I'g\|^2_{\dot H^1}e^{2\gamma \lambda^{-3}T'} +\frac{1}{\gamma^4}\|I'g\|^4_{L^2}e^{2\gamma \lambda^{-3}T'} \notag \\
  &  \ +|\int^{T'}_0(\Lambda_3+\Lambda_4)e^{2\gamma \lambda^{-3}t}dt| +|\int^{T'}_0\Lambda'_3e^{2\gamma \lambda^{-3}t}dt|^2), \label{lemma1 3}
\end{align}
where $\Lambda'_3= \Lambda_3((\xi_3+\eta_3)\frac{m(\zeta_1+\zeta_2)}{m(\zeta_1)m(\zeta_2)};I'v)$, $C_1>1$.
\end{lem}
{\bf Proof.}   Similarly to \eqref{E(I'v)exp}, one gets
\begin{align}
  \frac{d}{dt}\|I'v(t)\|^2_{L^2(\Real^2)}e^{2\gamma \lambda^{-3}t}& =2 \lambda^{-3}e^{2\gamma \lambda^{-3}t}\int_{\Real^2}I'gI'vdxdy -2e^{2\gamma \lambda^{-3}t}\Lambda'_3\nonumber
\end{align}
which implies
\begin{align}
  &\|I'v(t)\|^2_{L^2(\Real^2)}e^{2\gamma \lambda^{-3}t}\notag \\
\lesssim&  \|I'v_0\|^2_{L^2}+|\int^{T'}_0e^{2\gamma \lambda^{-3}t}\Lambda'_3dt|+\frac{e^{\gamma \lambda^{-3}T'}}{\gamma}\|I'g\|_{L^2} \mathop{sup}\limits_{t\in[0,T']}\|I'v(t)\|_{L^2}e^{\gamma \lambda^{-3}t} \notag \\
\lesssim&  \|I'v_0\|^2_{L^2}+|\int^{T'}_0e^{2\gamma \lambda^{-3}t}\Lambda'_3dt|+\frac{C_{\epsilon}}{\gamma^2}e^{2\gamma \lambda^{-3}T'}\|I'g\|^2_{L^2} +\epsilon \mathop{sup}\limits_{t\in[0,T']}\|I'v(t)\|^2_{L^2} e^{2\gamma \lambda^{-3}t}\label{lemma1 4}
\end{align}
for $t\in [0,T']$.

As the last term of \eqref{lemma1 4} can be absorbed by the left side via taking $\epsilon$ sufficiently small, we obtain \eqref{lemma1 2}.

Using H\"{o}rder's inequalities, one has
    \begin{align}
  &|\int^{T'}_0\int_{\Real^2}\lambda^{-3}I'g((\Delta-\partial_{x}\partial_{y})I'v+(I'v)^2 )e^{2\gamma \lambda^{-3}t}dxdydt| \notag \\
      \lesssim  &\frac{e^{\gamma \lambda^{-3}T'}}{\gamma}\|\nabla I'g\|_{L^2}\mathop{sup}\limits_{t\in[0,T']}\|\nabla I'v(t)\|_{L^2}e^{\gamma \lambda^{-3}t}
                  \notag \\
&\ +\frac{e^{\frac{2}{3}\gamma \lambda^{-3}T'}}{\gamma}\| I'g\|_{L^3}\mathop{sup}\limits_{t\in[0,T']}\|I'v(t)\|^2_{L^3}e^{\frac{4}{3}\gamma \lambda^{-3}t}\notag \\
      \lesssim & \frac{C_{\epsilon}e^{2\gamma \lambda^{-3}T'}}{\gamma^2}\|I'g\|^2_{\dot H^1}+\epsilon\mathop{sup}\limits_{t\in[0,T']}\|I'v(t)\|^2_{\dot H^1} e^{2\gamma \lambda^{-3}t}
\notag \\
&\  + \frac{e^{2\gamma \lambda^{-3}T'}}{\gamma^3}\| I'g\|^3_{L^3}+\mathop{sup}\limits_{t\in[0,T']}\|I'v(t)\|^3_{L^3}e^{2\gamma \lambda^{-3}t}. \label{Sec4 1}
    \end{align}
Hence, from \eqref{E(I'v)} and \eqref{Sec4 1} we obtain
\begin{align}
\|I'v(t)\|^2_{\dot H^1}e^{2\gamma \lambda^{-3}t}&\lesssim  \|I'v_0\|^2_{\dot H^1}+\frac{e^{2\gamma \lambda^{-3}T'}}{\gamma^2}\|I'g\|^2_{\dot H^1} +|\int^{T'}_0(\Lambda_3+\Lambda_4)e^{2\gamma \lambda^{-3}t}dt| \notag \\
      & \ +\frac{e^{2\gamma \lambda^{-3}T'}}{\gamma^3}\| I'g\|^3_{L^3}+ \mathop{sup}\limits_{t\in[0,T']}\|I'v(t)\|^3_{L^3}e^{2\gamma \lambda^{-3}t}  \label{Sec4 2}
\end{align}
for $t\in [0,T']$.

Gagliardo-Nirenberg's inequality and \eqref{lemma1 2} give
    \begin{align}
  \| I'g\|^3_{L^3} &\lesssim \| I'g\|_{\dot H^1}\| I'g\|^2_{L^2}\notag \\
 &\lesssim \gamma \| I'g\|^2_{\dot H^1}+\frac{1}{\gamma}\| I'g\|^4_{L^2}\label{Sec4 3}
\end{align}
    and
        \begin{align}
 & \mathop{sup}\limits_{t\in[0,T']} \| I'v\|^3_{L^3}e^{2\gamma \lambda^{-3}t}\notag \\
\lesssim&  \mathop{sup}\limits_{t\in[0,T']}\| I'v\|_{\dot H^1}\| I'v\|^2_{L^2}e^{2\gamma \lambda^{-3}t}\notag \\
\lesssim& C_{\epsilon} \mathop{sup}\limits_{t\in[0,T']}\| I'v\|^4_{L^2}e^{2\gamma \lambda^{-3}t}+\epsilon\mathop{sup}\limits_{t\in[0,T']}\| I'v\|^2_{\dot H^1}e^{2\gamma \lambda^{-3}t}\notag \\
\lesssim&  C_{\epsilon}\|I'v_0\|^4_{L^2}+\frac{C_{\epsilon}}{\gamma^4}e^{2\gamma \lambda^{-3}T'}\|I'g\|^4_{L^2}+C_{\epsilon} |\int^{T'}_0e^{2\gamma \lambda^{-3}t}\Lambda'_3dt|^2\notag \\
& + \epsilon \mathop{sup}\limits_{t\in[0,T']}\| I'v\|^2_{\dot H^1}e^{2\gamma \lambda^{-3}t}. \label{Sec4 4}
\end{align} 
Collecting \eqref{Sec4 2}, \eqref{Sec4 3} and \eqref{Sec4 4}, we can  get \eqref{lemma1 3}.

The following a-priori estimate is impactful for controlling $\|u\|_{H^s}$.
\begin{prop} \label {Prop priori estimate1}
  Let $C_2\ll1$, $C_3>1$, $C_4\gg1$ and $T'>0$ be given. Assume $v$ is a solution of  \eqref{rdfZK} on $[0,T']$.
  If $\lambda^{3}\geq \gamma $, $(N')^{\frac{1}{2}-}\geq C_4\lambda^{2}T'$,
  \begin{align}
1\lesssim\lambda^{2}\|I'v_0\|^2_{L^2}, \ \ 1\lesssim\lambda^{2}\|I'g\|^2_{L^2}, \nonumber
   \end{align}
  $$\|I'v_0\|^2_{L^2}+\frac{1}{\gamma^2}\|I'g\|^2_{L^2}e^{2\gamma \lambda^{-3}T'}\leq C_2,$$
  $$\|I'v_0\|^2_{\dot H^1}+\frac{1}{\gamma^2} \|I'g\|^2_{\dot H^1}e^{2\gamma \lambda^{-3}T'} \leq \lambda^{-2}C_2,$$
  then we have
\begin{align}
 \|I'v(T')\|^2_{L^2}e^{2\gamma \lambda^{-3}T'}\leq  C_3(\|I'v_0\|^2_{L^2}+\frac{1}{\gamma^2}\|I'g\|^2_{L^2}e^{2\gamma \lambda^{-3}T'}), \nonumber
   \end{align}
  \begin{align}
    \|I'v(T')\|^2_{\dot H^1}e^{2\gamma \lambda^{-3}T'}\leq& C_3 (\|I'v_0\|^2_{\dot H^1}+\|I'v_0\|^4_{L^2}+\frac{1}{\gamma^2} \|I'g\|^2_{\dot H^1}e^{2\gamma \lambda^{-3}T'}\notag \\
 & +\frac{1}{\gamma^4}\|I'g\|^4_{L^2}e^{2\gamma \lambda^{-3}T'}). \nonumber
 \end{align}
\end{prop}
{\bf Proof.} Actually choosing $\lambda$ sufficiently large, from scaling we have
\begin{align}
  &\|I'v_0\|^2_{\dot H^1}+\|I'v_0\|^4_{L^2}+\frac{1}{\gamma^2} \|I'g\|^2_{\dot H^1}e^{2\gamma \lambda^{-3}T'}+\frac{1}{\gamma^4}\|I'g\|^4_{L^2}e^{2\gamma \lambda^{-3}T'}\notag \\
\ll&\|I'v_0\|^2_{L^2}+\frac{1}{\gamma^2}\|I'g\|^2_{L^2}e^{2\gamma \lambda^{-3}T'}. \nonumber
 \end{align}
Taking $\delta=\delta(\|I'v_0 \|_{H^1}, \lambda^{-3}\|I'g \|_{H^1},\gamma\lambda^{-3})$ as in Proposition \ref {Prop lwp dfZK} and $j\in \mathbb{N}$ satisfying $\delta j=T'$. For $0\leq k \leq j$, $k\in \mathbb{N}$, we prove
\begin{align}
  &\|I'v(k\delta)\|^2_{L^2}e^{2\gamma \lambda^{-3}k\delta}\notag \\
\leq& 2C_1(\|I'v_0\|^2_{L^2}+\frac{1}{\gamma^2}\|I'g\|^2_{L^2}e^{2\gamma \lambda^{-3}k\delta})\notag \\
  \leq& 2C_1C_2 \label{prop prio 1}
\end{align}
and
  \begin{align}
    &\|I'v(k\delta)\|^2_{\dot H^1}e^{2\gamma \lambda^{-3}k\delta}\notag \\
    \leq& 2C_1 (\|I'v_0\|^2_{\dot H^1}+\|I'v_0\|^4_{L^2}+\frac{1}{\gamma^2} \|I'g\|^2_{\dot H^1}e^{2\gamma \lambda^{-3}k\delta}\notag \\
          & +\frac{1}{\gamma^4}\|I'g\|^4_{L^2}e^{2\gamma \lambda^{-3}k\delta})\notag \\
     \leq& 4C_1C_2 \label{prop prio 2}
  \end{align}
  by induction.

  For $k=0$, \eqref{prop prio 1} and \eqref{prop prio 2} hold trivially. We assume \eqref{prop prio 1} and \eqref{prop prio 2} hold true for $k=l$ where $0\leq l\leq j-1$. 
From Lemma \ref{I'v lemma1}, one has
\begin{align}
  &\|I'v((l+1)\delta)\|^2_{L^2}e^{2\gamma \lambda^{-3}(l+1)\delta}\notag \\
\leq&  C_1(\|I'v_0\|^2_{L^2}+\frac{1}{\gamma^2}\|I'g\|^2_{L^2}e^{2\gamma \lambda^{-3}(l+1)\delta}\ +|\int^{(l+1)\delta}_0e^{2\gamma \lambda^{-3}t}\Lambda'_3dt|) \nonumber 
\end{align}
and
\begin{align}
 &\|I'v((l+1)\delta)\|^2_{\dot H^1}e^{2\gamma \lambda^{-3}(l+1)\delta}\notag \\
\leq&  C_1 (\|I'v_0\|^2_{\dot H^1}+\|I'v_0\|^4_{L^2}+\frac{1}{\gamma^2} \|I'g\|^2_{\dot H^1}e^{2\gamma \lambda^{-3}(l+1)\delta} \notag \\
&  \  +\frac{1}{\gamma^4}\|I'g\|^4_{L^2}e^{2\gamma \lambda^{-3}(l+1)\delta}+|\int^{(l+1)\delta}_0(\Lambda_3+\Lambda_4)e^{2\gamma \lambda^{-3}t}dt|\notag \\
&  \  +|\int^{(l+1)\delta}_0\Lambda'_3e^{2\gamma \lambda^{-3}t}dt|^2).\nonumber
\end{align}
Therefore, it suffices to prove
\begin{align}
|\int^{(l+1)\delta}_0e^{2\gamma \lambda^{-3}t}\Lambda'_3dt|\leq  \|I'v_0\|^2_{L^2}+\frac{1}{\gamma^2}\|I'g\|^2_{L^2}e^{2\gamma \lambda^{-3}(l+1)\delta} \label{Lambda 3a}
\end{align}
and
\begin{align}
  &|\int^{(l+1)\delta}_0(\Lambda_3+\Lambda_4)e^{2\gamma \lambda^{-3}t}dt| +|\int^{(l+1)\delta}_0e^{2\gamma \lambda^{-3}t}\Lambda'_3dt|^2\notag \\
  \leq& \|I'v_0\|^2_{\dot H^1}+\|I'v_0\|^4_{L^2}+\frac{1}{\gamma^2} \|I'g\|^2_{\dot H^1}e^{2\gamma \lambda^{-3}(l+1)\delta} +\frac{1}{\gamma^4}\|I'g\|^4_{L^2}e^{2\gamma \lambda^{-3}(l+1)\delta}.\label{Lambda 3and4}
\end{align}
From Proposition \ref{p5}, we have
\begin{align}
&|\int^{(l+1)\delta}_0e^{2\gamma \lambda^{-3}t}\Lambda'_3dt|\notag \\
  =&|\int^{(l+1)\delta}_0 \Lambda_3((\xi_3+\eta_3)\frac{m(\zeta_1+\zeta_2)}{m(\zeta_1)m(\zeta_2)};e^{\frac{2}{3}\gamma \lambda^{-3}t}I'v)|\notag \\
  \lesssim& \delta^{\frac{2}{3}}(N')^{-2+}\sum^l_{m=0}\|e^{\frac{2}{3}\gamma \lambda^{-3}t}I'v\|^3_{Y_{[m\delta,(m+1)\delta]}^1}\notag \\
   \lesssim& \delta^{\frac{2}{3}}(T')^{-1}\sum^l_{m=0}\|I'v\|^3_{Y_{[m\delta,(m+1)\delta]}^1}e^{2\gamma \lambda^{-3}(m+1)\delta}   \label{Lambda(3a1)}.
\end{align}
Proposition \ref {Prop lwp dfZK}, \eqref{prop prio 1} and \eqref{prop prio 2} provide the up bound of $\|I'v\|_{Y^1}$
\begin{align}
  &\|I'v\|^3_{Y_{[m\delta,(m+1)\delta]}^1}e^{2\gamma \lambda^{-3}(m+1)\delta}\notag \\
  \lesssim &\|I'v(m\delta)\|^3_{H^1}e^{2\gamma \lambda^{-3}(m+1)\delta}+(\lambda^{-3}\|I'g\|_{H^1})^3e^{2\gamma \lambda^{-3}(m+1)\delta}\notag \\
 \lesssim &(C_1C_2)^{\frac{1}{2}}C_1(\|I'v_0\|^2_{L^2}+\frac{1}{\gamma^2}\|I'g\|^2_{L^2}e^{2\gamma \lambda^{-3}(m+1)\delta}) \label{Lambda(3a2)}.
\end{align}
Hence, from \eqref{Lambda(3a1)} and \eqref{Lambda(3a2)}, one obtains
\begin{align}
&|\int^{(l+1)\delta}_0e^{2\gamma \lambda^{-3}t}\Lambda'_3dt|\notag \\
  \lesssim& (C_1C_2)^{\frac{1}{2}}C_1\delta^{\frac{2}{3}}(T')^{-1}(l+1)(\|I'v_0\|^2_{L^2}+\frac{1}{\gamma^2}\|I'g\|^2_{L^2}e^{2\gamma \lambda^{-3}(l+1)\delta})\notag \\
\lesssim& (C_1C_2)^{\frac{1}{2}}C_1\delta^{-\frac{1}{3}}(\|I'v_0\|^2_{L^2}+\frac{1}{\gamma^2}\|I'g\|^2_{L^2}e^{2\gamma \lambda^{-3}(l+1)\delta}) \label{Lambda(3a3)}.
\end{align}
We can get \eqref{Lambda 3a} from \eqref{Lambda(3a3)} as $C_2\ll1$.

Moreover,
\begin{align}
&|\int^{(l+1)\delta}_0e^{2\gamma \lambda^{-3}t}\Lambda'_3dt|^2\notag \\
  \leq& \frac{1}{4}(\|I'v_0\|^2_{\dot H^1}+\|I'v_0\|^4_{L^2}+\frac{1}{\gamma^2} \|I'g\|^2_{\dot H^1}e^{2\gamma \lambda^{-3}(l+1)\delta} +\frac{1}{\gamma^4}\|I'g\|^4_{L^2}e^{2\gamma \lambda^{-3}(l+1)\delta})\nonumber
\end{align}
which allows that one only need to show
\begin{align}
  &|\int^{(l+1)\delta}_0(\Lambda_3+\Lambda_4)e^{2\gamma \lambda^{-3}t}dt| \notag \\
  \leq& \frac{1}{2}(\|I'v_0\|^2_{\dot H^1}+\|I'v_0\|^4_{L^2}+\frac{1}{\gamma^2} \|I'g\|^2_{\dot H^1}e^{2\gamma \lambda^{-3}(l+1)\delta} +\frac{1}{\gamma^4}\|I'g\|^4_{L^2}e^{2\gamma \lambda^{-3}(l+1)\delta}).\label{Lambda 3and4 1}
\end{align}
Similarly, from Proposition \ref{p3} and \ref{p4}, we have
\begin{align}
  &|\int^{(l+1)\delta}_0(\Lambda_3+\Lambda_4)e^{2\gamma \lambda^{-3}t}dt| \notag \\
  \lesssim& \delta^{\frac{1}{6}}(N')^{-\frac{1}{2}+}\sum^l_{m=0}\|I'v\|^3_{Y_{[m\delta,(m+1)\delta]}^1}e^{2\gamma \lambda^{-3}(m+1)\delta} \notag \\
  \lesssim& \delta^{\frac{1}{6}}(C_4\lambda^2T')^{-1}\sum^l_{m=0}(\|I'v(m\delta)\|^3_{H^1}+(\lambda^{-3}\|I'g\|_{H^1})^3)e^{2\gamma \lambda^{-3}(m+1)\delta} .\label{Lambda 3and4 2}
\end{align}

On one hand, \eqref{prop prio 1} and \eqref{prop prio 2} tell us that
\begin{align}
  &\|I'v(m\delta)\|^3_{L^2}e^{2\gamma \lambda^{-3}(m+1)\delta}\notag \\
\lesssim& C^2_1(\|I'v_0\|^3_{L^2}+\frac{1}{\gamma^2}\|I'g\|^3_{L^2}e^{2\gamma \lambda^{-3}(m+1)\delta})\notag \\
\lesssim&  C^2_1\lambda^2(\|I'v_0\|^4_{L^2}+\frac{1}{\gamma^4}\|I'g\|^4_{L^2}e^{2\gamma \lambda^{-3}(m+1)\delta})\label{Lambda 3and4 3}
\end{align}
and
\begin{align}
  &\|I'v(m\delta)\|^3_{\dot H^1}e^{2\gamma \lambda^{-3}(m+1)\delta}\notag \\
  \lesssim& (C_1C_2)^{\frac{1}{2}} (\|I'v_0\|^2_{\dot H^1}+\|I'v_0\|^4_{L^2}+\frac{1}{\gamma^2} \|I'g\|^2_{\dot H^1}e^{2\gamma \lambda^{-3}(m+1)\delta}\notag \\
&  +\frac{1}{\gamma^4}\|I'g\|^4_{L^2}e^{2\gamma \lambda^{-3}(m+1)\delta}), \label{Lambda 3and4 4}
  \end{align}
on the other hand, we have
  \begin{align}
  &(\lambda^{-3}\|I'g\|_{H^1})^3e^{2\gamma \lambda^{-3}(m+1)\delta}\notag \\
  \lesssim& \frac{1}{\gamma^3}\|I'g\|^3_{L^2}e^{2\gamma \lambda^{-3}(m+1)\delta}+ \frac{1}{\gamma^3}\|I'g\|^3_{\dot H^1}e^{2\gamma \lambda^{-3}(m+1)\delta}\notag \\
    \lesssim& \lambda^2\frac{1}{\gamma^4}\|I'g\|^4_{L^2}e^{2\gamma \lambda^{-3}(m+1)\delta}+ (C_2)^{\frac{1}{2}} \frac{1}{\gamma^2}\|I'g\|^2_{\dot H^1}e^{2\gamma \lambda^{-3}(m+1)\delta}\notag \\
   \lesssim&  \lambda^2(\frac{1}{\gamma^4}\|I'g\|^4_{L^2}e^{2\gamma \lambda^{-3}(m+1)\delta}+  \frac{1}{\gamma^2}\|I'g\|^2_{\dot H^1}e^{2\gamma \lambda^{-3}(m+1)\delta}). \label{Lambda 3and4 5}
  \end{align}
  Collecting \eqref{Lambda 3and4 2}-\eqref{Lambda 3and4 5}, one gets
  \begin{align}
  &|\int^{(l+1)\delta}_0(\Lambda_3+\Lambda_4)e^{2\gamma \lambda^{-3}t}dt| \notag \\
    \lesssim& \delta^{\frac{1}{6}}(C_4\lambda^2T')^{-1} C^2_1\lambda^2(l+1)(\|I'v_0\|^2_{\dot H^1}+\|I'v_0\|^4_{L^2}\notag \\
  & +\frac{1}{\gamma^2} \|I'g\|^2_{\dot H^1}e^{2\gamma \lambda^{-3}(l+1)\delta} +\frac{1}{\gamma^4}\|I'g\|^4_{L^2}e^{2\gamma \lambda^{-3}(l+1)\delta})\notag \\
 \lesssim& \delta^{-\frac{5}{6}}C^2_1C^{-1}_4 (\|I'v_0\|^2_{\dot H^1}+\|I'v_0\|^4_{L^2}\notag \\
  & +\frac{1}{\gamma^2} \|I'g\|^2_{\dot H^1}e^{2\gamma \lambda^{-3}(l+1)\delta} +\frac{1}{\gamma^4}\|I'g\|^4_{L^2}e^{2\gamma \lambda^{-3}(l+1)\delta})  \nonumber
  \end{align}
  which gives \eqref{Lambda 3and4 1} by taking $C_4$ sufficiently large.  
\begin{prop} \label {Prop priori estimate2}
Let $C_2\ll1$, $C_3>1$, $C_4\gg1$ and $T>0$ be given. Assume $u$ is a solution of  \eqref{dfZK} on $[0,T]$.
  If $N^{\frac{3}{11}-}\geq\gamma $, $N^{0+}\geq C_4T$, 
  $$\|Iu_0\|^2_{L^2}+\frac{1}{\gamma^2}\|If\|^2_{L^2}e^{2\gamma T}\leq C_2N^{\frac{2}{11}-},$$
  and
  $$\|Iu_0\|^2_{\dot H^1}+\frac{1}{\gamma^2} \|If\|^2_{\dot H^1}e^{2\gamma T} \leq C_2N^{\frac{2}{11}-},$$
  then we have
\begin{align}
 \|Iu(T)\|^2_{L^2}e^{2\gamma T}\leq  C_3(\|Iu_0\|^2_{L^2}+\frac{1}{\gamma^2}\|If\|^2_{L^2}e^{2\gamma T}), \nonumber
   \end{align} 
  \begin{align}
\|Iu(T)\|^2_{\dot H^1}e^{2\gamma T}&\leq C_3 (\|Iu_0\|^2_{\dot H^1}+\|Iu_0\|^4_{L^2}+\frac{1}{\gamma^2} \|If\|^2_{\dot H^1}e^{2\gamma T}+\frac{1}{\gamma^4}\|If\|^4_{L^2}e^{2\gamma T}).\nonumber
 \end{align}
\end{prop}
{\bf Proof.} Notice that $\|I'v\|^2_{L^2}= \lambda^{-2}\|Iu\|^2_{L^2}$, $|I'g\|^2_{L^2}= \lambda^{-2}\|If\|^2_{L^2}$, $\|I'v\|^2_{\dot H^1}= \lambda^{-4}\|Iu\|^2_{\dot H^1}$, $\|I'g\|^2_{\dot H^1}= \lambda^{-4}\|If\|^2_{\dot H^1}$ and the fact $\lambda^{2}\|I'v_0\|^2_{L^2}= \|Iu_0\|^2_{L^2}>\frac{1}{2}\|u_0\|^2_{L^2}$ if $N$ is sufficiently large, we acquire Proposition \ref{Prop priori estimate2} via taking $\lambda=N^{\frac{1}{11}-}$, $N'=\frac{N}{\lambda}$ and $T'=\lambda^{3}T$ in Proposition \ref{Prop priori estimate1}. 

Finally ,we show the exsitence of the global attractor.

{\bf Proof of Theorem \ref{Global Attactor}.} 
Let $0 < \epsilon \ll 1$ be fixed. We choose  $T_{1} > 0$ so that 
\begin{align} 
  e^{2 \gamma T_{1}} > &(\|u_{0}\|^{2}_{H^{s}}+\|u_0\|^4_{L^2})(\frac{1}{\gamma^2} \|f\|^2_{ H^1}+\frac{1}{\gamma^4}\|f\|^4_{L^2})^{-1} \mathop{max} \bigg\lbrace \gamma^{\frac{22(1-s)}{3}+}, (C_{4} T_{1})^{\frac{2(1-s)}{\epsilon}},\notag \\
                       &  \left(\frac{C_{2}}{2}\|u_{0}\|^{-2}_{H^{s}}\right)^{\frac{11(1-s)}{(10-11s)}-} , \left( 2 C_{2}^{-1} \gamma^{-2} \|f\|^{2}_{H^{1}}e^{2\gamma T_{1}} \right)^{11(1-s)+} \bigg\rbrace , \label{TH1.2 1}
\end{align}
which is possible as $11(1-s)<1$.  $T_{1}$ depends only on $\|u_{0}\|_{H^{s}}$, $\|f\|_{H^{1}}$ and $\gamma$. 
Set
\begin{align} \label{TH1.2  2}
N = \mathop{max} \bigg\lbrace \gamma^{\frac{11}{3}+},(C_{4} T_{1})^{\frac{1}{\epsilon}}, \left(  \frac{C_{2}}{2}\|u_{0}\|^{-2}_{H^{s}}\right)^{\frac{11}{2(10-11s)}-} , 
\left( 2 C_{2}^{-1} \gamma^{-2} \|f\|^{2}_{H^{1}} e^{2\gamma T_{1}} \right)^{\frac{11}{2}+} \bigg\rbrace.
\end{align}
From the choice of $T_1$ and $N$, we know
$$N^{\frac{3}{11}-} \geq \gamma, \ \  \ \ \ \ N^{\epsilon} \geq C_{4}T_{1},$$
and 
$$\|Iu_{0}\|_{H^{1}}^{2} \leq N^{2-2s} \|u_{0}\|_{H^{s}}^{2}  \leq \dfrac{C_{2}}{2} N^{\frac{2}{11}-},$$
$$\gamma^{-2} \|If\|^{2}_{H^{1}} e^{2\gamma T_{1}}\leq \dfrac{C_{2}}{2} N^{\frac{2}{11}-}.$$
Hence, from Proposition \ref{Prop priori estimate2}, one gains
\begin{align}
  \|u(T_{1})\|_{H^{s}}^{2} \leq& \|Iu(T_{1})\|_{H^{1}}^{2}\notag \\
\leq &C_3 (\|Iu_0\|^2_{ H^1}e^{-2\gamma T}+\|Iu_0\|^4_{L^2}e^{-2\gamma T}+\frac{1}{\gamma^2} \|If\|^2_{ H^1}+\frac{1}{\gamma^4}\|If\|^4_{L^2})\notag \\
\leq &C_3 (N^{2(1-s)}\|u_0\|^2_{ H^s}e^{-2\gamma T}+\|u_0\|^4_{L^2}e^{-2\gamma T}+\frac{1}{\gamma^2} \|f\|^2_{ H^1}+\frac{1}{\gamma^4}\|f\|^4_{L^2}).\nonumber                                   \end{align}
From \eqref{TH1.2 1} and \eqref{TH1.2 2} , we get
$$N^{2(1-s)}e^{-2\gamma T_1} (\|u_0\|^2_{ H^s}+\|u_0\|^4_{L^2})< \frac{1}{\gamma^2} \|f\|^2_{ H^1}+\frac{1}{\gamma^4}\|f\|^4_{L^2}$$
which helps us give the bound
$$\|u(T_{1})\|_{H^{s}}^{2} \leq 2C_{3}( \frac{1}{\gamma^2} \|f\|^2_{ H^1}+\frac{1}{\gamma^4}\|f\|^4_{L^2})< K_{1}$$
where $K_{1}$ depends only on $\|f\|_{ H^1}$ and $\gamma$.

In the next place, one can fix $T_{2} > 0$ and solve  \eqref{dfZK} on time interval $[T_{1},T_{1} + T_{2}]$ with initial data replaced by $u(T_{1}).$ Let $K_{2} > 0$ be sufficiently large such that
\begin{align} 
  K_2e^{2 \gamma t} > & \mathop{max} \bigg\lbrace \gamma^{\frac{22(1-s)}{3}+}, (C_{4} t)^{\frac{2(1-s)}{\epsilon}}, \left(2C_2^{-1}K_1\right)^{\frac{11(1-s)}{(11s-10)}+} ,\notag \\
  & \ \ \ \left( 2 C_{2}^{-1} \gamma^{-2} \|f\|^{2}_{H^{1}}e^{2\gamma t} \right)^{11(1-s)+} \bigg\rbrace , \label{TH1.2 3}
\end{align}
for any $t > 0$. Set $N^{2(1-s)} = K_{2}e^{2\gamma T_{2}}$, then inequality \eqref{TH1.2 3} verifies the assumptions in Proposition \ref{Prop priori estimate2} and hence we obtain
\begin{align}
\|Iu(T_{1} + T_{2})\|_{H^{1}}^{2} \leq &C_3 (N^{2(1-s)}\|u(T_1)\|^2_{ H^s}e^{-2\gamma T_2}+\|u(T_1)\|^4_{L^2}e^{-2\gamma T_2}+\frac{1}{\gamma^2} \|f\|^2_{ H^1}+\frac{1}{\gamma^4}\|f\|^4_{L^2}) \notag \\
 \leq &C_3 (K_1K_2+K^2_1+\frac{1}{\gamma^2} \|f\|^2_{ H^1}+\frac{1}{\gamma^4}\|f\|^4_{L^2})<K_3. \nonumber
\end{align}

For $t > T_{1},$ we define the maps $M_{1}(t)$ and $M_{2}(t)$ as
$$\widehat{M_{1}(t)u_{0}} = \widehat{A(t)u_{0}}|_{|\zeta| < N_t}, \ \ \  \widehat{M_{2}(t)u_{0}} = \widehat{A(t)u_{0}}|_{|\zeta| > N_t},$$
where $A(t)u_{0} = u(t)$ and $N_t = (K_{2}e^{2\gamma(t-T_{1})})^{\frac{1}{2(1-s)}}.$

It's easy to see that for $t > T_{1},$ 
\begin{align}
\|M_{1}(t)u_{0}\|_{H^{1}}^{2} \leq \|Iu(t)\|_{H^{1}}^{2} &< K_{3}, \notag \\
\|M_{2}(t)u_{0}\|_{H^{s}}^{2} \leq N^{2s-2}\|Iu(t)\|_{H^{1}}^{2} &< K_{2}^{-1} K_{3}e^{-2\gamma(t - T_{1})}.\nonumber
\end{align}
 Hence we obtain Theorem \ref{Global Attactor} by taking $K = \mathop{max}\lbrace K_{3}^{\frac{1}{2}}, K_{2}^{-\frac{1}{2}}K_{3}^{\frac{1}{2}}\rbrace$.
\vspace{4mm}
\\{\bf Acknowledgment}\vspace{2mm}
\\The author would like to express his deep gratitude to Professor Yoshio Tsutsumi and Professor Baoxiang Wang for their contructive guidance and kind help.

\bibliographystyle{amsplain}
\bibliography{}

\end{document}